            \documentstyle[12pt]{amsart}
                        \numberwithin{equation}{section}
                        \theoremstyle{plain}
                         \newtheorem{thm}{Theorem}[section]
                          \newtheorem{lem}[thm]{Lemma}
                                \newtheorem{pro}[thm]{Proposition}
                                \newtheorem{cor}[thm]{Corollary}
                       \newtheorem{conj}[thm]{Conjecture}
                         \theoremstyle{definition}
                                \newtheorem{rem}{Remark}[section]

\begin{document}

                  \title[On perturbative $PSU(N)$ invariants]
{On perturbative $PSU(N)$ invariants of 
rational homology 3-spheres}

                \author[ Thang Le ]{Thang T. Q. Le }
\thanks{This work is partially supported by NSF grant DMS-9626404.
This and related preprints can be obtained at 
{\tt  http://\linebreak[0]\linebreak[0]www.\linebreak[0]math.\linebreak
[0]buffalo.\linebreak[0]edu/\linebreak[0]$\sim$letu/}}
\address{Dept. of Mathematics, SUNY at Buffalo, Buffalo, NY 14214, USA}
\email{letu@@newton.math.buffalo.edu}

\begin{abstract} 
We construct power series invariants of rational
homology 3-spheres from  quantum $PSU(n)$-invariants. 
The power series can be regarded as perturbative invariants
corresponding to the contribution of the
trivial connection in the hypothetical Witten's integral.
This generalizes a result of Ohtsuki (the $n=2$ case) which led him
to the definition of finite type invariants of 3-manifolds. The proof utilizes
some symmetry properties of quantum invariants (of links) derived from
the theory of affine Lie algebras and the theory of the Kontsevich integral.    
                  \end{abstract}
\maketitle

\newcommand{\ab}{{\bold a}}
\newcommand{\boldeta}{{\boldsymbol \eta}}

\addtocounter{section}{-1}
\newcommand{\gau}{\gamma}
\newcommand{\Jhom}{J^{Hom}}
\newcommand{\bt}{{\bold t}}
\newcommand{\bb}{{\bold b}}
\newcommand{\tJ}{\tilde J}
\newcommand{\Lroot}{\Lambda^{\text{root}}}\newcommand{\RRR}{\Lroot}
\newcommand{\bk}{{\bold k}}
\newcommand{\bs}{{\bold s}}
\newcommand{\arcsinh}{\operatorname{arcsinh}}
\newcommand{\taut}{\tau}
\newcommand{\tauone}{\tau^{(1)}}
\newcommand{\tautwo}{\tau^{(2)}}
\newcommand{\tauthree}{\tau^{(3)}}
\newcommand{\tauPSU}{\tau^{PSU(n)}}
\newcommand{\Weyl}{{\cal W}}
\newcommand{\Qtwo}{Q_L^{(2)}}
\newcommand{\Qthree}{Q_L^{(3)}}
\newcommand{\Weylr}{\Weyl_{(r)}}
\newcommand{\cF}{{\cal F}}
\newcommand{\Zplus}{{\Bbb Z}_+}
 \newcommand{\Jo}{J^0} \newcommand{\GG}{{\cal G}}
\newcommand{\tF}{\tilde F}\newcommand{\Qo}{{\overset{o}{Q}}}
\newcommand{\Zr}{{\Bbb Z}_{(r)}}
\newcommand{\Root}{R_+}
\newcommand{\sn}{\operatorname{sn}}
\newcommand{\SSS}{{C_r}}
\newcommand{\RR}{{\cal R}}
\newcommand{\Po}{\Lambda_{+}} \newcommand{\Pnr}{P(N,r-N)} 
\newcommand{\tPnr}{\tilde\Pnr}
\newcommand{\QI}{\Qo^{(1)}}
 \newcommand{\QII}{\Qo^{(2)}}
\newcommand{\QIII}{\Qo^{(3)}}            
\newcommand{\G}{\Gamma}
\newcommand{\Q}{{\Bbb Q}}
\newcommand{\suml} {{\sum{}^{(\ell)}}}
\newcommand{\cAm}{{\cal A}(\sqcup^m S^1)}
\newcommand{\Zell}{{\Bbb Z}_{(\ell)}}
\newcommand{\tQ}{\tilde Q}
\newcommand{\bR}{{\Bbb R}}
\newcommand{\CP}{{\cal Q}}
\newcommand{\fS}{\SS}
\newcommand{\Sa}{{\cal S}}

\newcommand{\cP}{{\cal P}}
\newcommand{\cC}{{\cal B}}
\newcommand{\cH}{{\cal H}} \newcommand{\al}{\alpha}
\newcommand{\Qone}{Q_L^{(1)}}
\newcommand{\Fone}{F_{L'}^{(1)}}
\newcommand{\Poo}{\Lambda_{++}}
\newcommand{\CC}{{\cal F}}
\newcommand{\ZCC}{{\cal F_{\Bbb Z}}}
\newcommand{\ZG}{{\cal G_{\Bbb Z}}}
\newcommand{\cA}{{\cal A}}
\newcommand{\bl}{{\bold l}}

\newcommand{\hZ}{{\Hat Z}}\newcommand{\cZ}{{\check Z}}
\newcommand{\ZA}{{\cA^{\Bbb Z}}}
\newcommand{\ZP}{{\cP^{\Bbb Z}}}
\newcommand{\R}{{\Bbb R}}\newcommand{\D}{{\cal D}}
\newcommand{\ve}{\varepsilon}
\newcommand{\Z}{{\Bbb Z}}\newcommand{\C}{{\Bbb C}}
\newcommand{\cB}{{\cal B}}\newcommand{\Do}{\overset {\circ}{\cal D}}
\newcommand{\kk}{{\bold k}}
\newcommand{\Lrootr}{\Lambda^{\text{root}}_r}
\newcommand{\summu}{\sum_{\mu\in(\rho+\Lrootr)}}
\newcommand{\summuj}{\sum_{\mu_j\in(\rho+\Lrootr)}}
\newcommand{\Ftwo}{F_{L'}^{(2)}}

\newcommand{\Fthree}{F_{L'}^{(3)}}
\newcommand{\Rin}{\Q[[x]]_{(M)}}

        \section{Introduction}
In this paper we construct  power series invariants  of rational
homology 3-spheres from the quantum $PSU(n)$-invariants . 
This generalizes a result of Ohtsuki (the $n=2$ case) which led him
to the definition of finite type invariants of 3-manifolds.

For a fixed compact Lie group and an integer $r$, 
Witten defined an invariant of 3-manifolds using Feynman path integral which
is not mathematically rigorous. Reshetikhin and Turaev  \cite{RT2}, and later
Turaev and Wenzl \cite{TuraevWenzl} defined (mathematically) quantum 
invariants of 3-manifolds associated
with simple Lie algebras (see \cite{Turaev} and
references therein). The quantum invariants can be thought of
as a mathematical realization of Witten's path integral.

One could  use formal perturbation theory to 
approximate the Witten integral. The coefficients of the approximations are 
known as perturbative invariants, which also need a  mathematically 
rigorous definition. 
One might want to  get perturbative invariants directly 
from (Reshetikhin-Turaev)
 quantum invariants.
The difficulty here is that quantum invariant can be defined only when 
the ``quantum'' parameter $q$ is  roots of
unity, while in order to get perturbative invariants, one need to make the formal
expansion $q= \exp(h)$, with $h$ an {\em indeterminate}.
 This substitution does not make sense if
$q$ takes only values roots of unity.

Ohtsuki \cite{Ohtsuki1} showed that one can   extract
 power series invariants (of  homology 3-spheres)
from quantum $SU(2)$-invariants which can be considered as the perturbative invariants. 
Ohtsuki's result led him to the important definition of finite type invariant of
3-manifolds, a counterpart of Vassiliev invariants for homology 3-spheres. The Ohtsuki series
can be considered, on the physics level, as the contribution of the trivial connection
 in the perturbative
expansion of the Witten integral (see \cite{Rozansky}). Further investigation of
the Ohtsuki's series was carried out by Ohtsuki \cite{Ohtsuki1,Ohtsuki3},
Lawrence and  Rozansky \cite{Lawrence,LawrenceRozansky,Rozansky,Rozansky2}, 
Lin and Wang \cite{LinWang},
Kricker and Spence \cite{Kricker}, and others. A result of H. Murakami 
\cite{Murakami} says that the first
term of the Ohtsuki series is 6 times the Casson-Walker invariant.

In \cite{LMO} it was conjectured that similar power series invariants exist
for Lie group $SU(n), n>2$. It is the main goal of this paper to prove this conjecture.
We show that one can extract power series invariants of {\em rational homology
3-spheres} from quantum $PSU(n)$-invariants. Here the quantum $PSU(n)$-invariant
is a refined version of the $SU(n)$-invariant which was introduced
by Kirby and Melvin for the $n=2$ case (see \cite{KirbyMelvin}), and by Kohno and Takata 
for the $n>2$ case (see \cite{KohnoTakata2}).
It will be clear from the proof that for rational homology 3-spheres, 
one should  use the $PSU(n)$ version instead of the $SU(n)$ one. It is because of some 
integrality and symmetry properties of quantum link invariants related to the $PSU(n)$ case.
In a sense, the $SU(n)$ version is based on the weight lattice, and the $PSU(n)$ version
is based on the root lattice; and quantum invariants of links on the root lattice are
nicer: they have integrality property and more symmetry.

A few words about the proof.
The proof of the main result differs from that of Ohtsuki for the case $n=2$, since
Ohtsuki \cite{Ohtsuki1} used some identities which are either specific to the $n=2$
case or hard to generalize to the $n>2$ case
(see also \cite{LinWang,Rozansky}). Some of our arguments are generalizations
of those in \cite{Ohtsuki1,LinWang,Rozansky} which treat the $n=2$ case.
We will make use of some symmetry properties
of quantum invariants of links which are derived from affine Lie algebra theory and
the Kontsevich integral theory. Along the way we also investigate the dependence of
quantum invariants of links on the colors of the link components. The proof also suggests
how to generalize the result to other Lie algebras, and we will try to formulate
everything in terms of Lie algebra theory.

In \S 1 we introduce the notation and formulate the main result. In section 2
we give the proof of the main theorem using some results proved in sections 3 and 4.
In section 3 we investigate the dependence of quantum invariants of links on the colors
(which are modules of Lie algebras). 
In section 4 we prove some properties of the quadratic Gauss sums
based on the root lattice.
In Appendix we show that our definition of quantum $PSU(n)$-invariant coincides
with that of Kohno and Takata, and calculate the value of a multi-variable
quadratic Gauss sum.

{\bf Acknowledgment.} Much of this work was carried out while the author was 
visiting  the Mathematical Sciences Research Institute in Berkeley in 1996-1997. 
Research at MSRI is supported in part by NSF grant DMS-9022140. The author thanks
J. Murakami, T. Takata and H. Wenzl for helpful  discussions. 

\section{Preliminaries. Main result}
\subsection{Notations}

We will use the following notations for objects related to the Lie algebra
$sl_n$. See, for example, \cite{Humphreys,Kac}, for the terminologies used  here.

$\al_1,\dots,\al_{n-1}$: \quad      the standard root basis.

$\Phi$: \quad the root system.

$\Phi_+$: \quad  the set of positive roots. There are $n(n-1)/2$ positive roots.

$\lambda_1,\dots,\lambda_{n-1}$: \quad the standard fundamental weights. 

$\Lambda$:  \quad the integral weight lattice, $\Lambda
=\{\,\sum_{i=1}^{n-1} k_i\lambda_i, k_i\in\Z\,\}$. All the weights $\lambda_i$ and 
roots $\alpha_j$ are  elements of $\Lambda$. In 
$\Lambda\otimes \R$ there is the  
standard scalar product, for which $(\alpha_i|\alpha_j)=A_{ij}$. Here $A$ is the 
Cartan matrix. One has $(\alpha_i|\lambda_j)= \delta_{ij}$.

$\rho$: \quad half-sum of positive roots. One has that
 $\rho= \lambda_1+\dots+\lambda_{n-1}$, 
and $(\rho|\rho)= n(n^2-1)/12$.

$\Zplus$: \quad the set of non-negative integers.

$\Po$: \quad the set of dominant weights, $\Po = \{\,k_1\lambda_1+\dots + 
k_{n-1}\lambda_{n-1}
 \mid  k_1, 
\dots,k_{n-1} \in\Zplus\,\}$

$\Poo$: \quad the set of strongly dominant weights, 
$\Poo = \rho+ \Po$.

$\Lroot$: \quad The root lattice, 
$\Lroot
=\{\,\sum_{i=1}^{n-1} k_i\al_i, k_i\in\Z\,\}$. The root lattice is a subgroup of
the weight lattice $\Lambda$ of index $n$.

$\theta$: \quad the longest root, $\theta = \alpha_1+\dots+\alpha_{n-1}.$

$\RR$:  \quad the ring of finite dimensional  $sl_n$-modules.

$V_\mu$, with $\mu\in\Poo$: \quad the irreducible $sl_n$-module with highest 
weight $\lambda-\rho$. Note that this is slightly different
from the usual notation. The shift by $\rho$ is more convenient for us.

The Weyl group $\Weyl$ acts on $\Lambda\otimes\R$, and a fundamental domain is the
fundamental chamber
$$C=\{\,z\in \Lambda\times\R \mid (z,\al_i)\ge 0,\quad i=1,\dots,n-1\,\} =
\{ \sum_{i=1}^{n-1} a_i \lambda_i \mid a_i \ge 0, a_i\in\R\}.$$
So, $\Po$ is the intersection of $\Lambda$ and $C$. The Weyl group
is generated by reflections along the boundary facets of $C$.
For an element $w\in\Weyl$, let $\sn(w)$ be the sign of $w$. Similarly,
for a non-zero number $b$, let $\sn(b)$ be the sign of $b$.

Throghout the paper, the number $n$ (of $sl_n$) is fixed. 
We will use $r$ to 
denote the ``level"  of
the theory (it's just a positive integer). Let
$$\bar r =  \frac{r-1 -n(n-1)}{2}.$$
 Note that if $r$ is odd,
which we always assume, then $\bar r$ is an integer.

We will use $q,x$ as indeterminates, with $q= 1+x$. We use $\zeta$ to denote
the $r$-th root of unity, $\zeta = \exp(2\pi i/r)$.
When $r$ is an odd prime, we will identify 
$\Z[\zeta]$ with
$\Z[q]/(q^{r-1}+\dots+q+1)$.

For any positive integer $\ell$ let 
$$\Zell =\Z[\frac{1}{(\ell-1)! \, n!}].$$ 
For a rational homology 3-sphere $M$ (i.e. a closed oriented 3-manifold
whose homology group $H_1(M,\Z)$ is finite), let
$\Rin$ be the set of all power series $\sum_{\ell} d_\ell x^\ell$ such that
$d_\ell$ is a number in $\Z[\frac{1}{(2\ell+ n(n-1))!\,|H_1(M,\Z)| }]$.

\subsection{Quantum invariants of oriented framed links and 3-manifolds}
For a framed oriented link $L$ with $m$ ordered components in $S^3$
 let $J_L(\mu_1,
\dots, \mu_m;q)$ be the $sl_n$-quantum invariant of the link (see
\cite{RT1,Turaev}), where
the components of $L$ are colored by $V_{\mu_1},\dots,V_{\mu_m}$.
Here  $q$ is an indeterminate and
the $\mu_j$'s are in $\Poo$, so that $V_{\mu_j}$ 
are finite-dimensional irreducible $sl_n$-modules of highest weight $\mu_j-\rho$.
Actually, $J_L$ can be regarded as a linear homomorphism from $\RR^
{\otimes m}$ to 
$\Z[q^{\pm 1/2n}]$.

We extend the definition of $J_L(\mu_1,\dots,\mu_m;q)$ to the case
when the $\mu_j$'s are 
in $\Lambda$ as follows. If one of the $\mu_j$ is on the boundary of
a chamber
$w(C)$, where $w$ is in the Weyl group, put $J_L(\mu_1,\dots,\mu_m;q)=0$. 
Otherwise, there exist
unique $w_1,\dots,w_m$ in the Weyl group $\Weyl$ such that the elements
$w_j(\mu_j)$ are in $\Poo$. Then put 
$$J_L(\mu_1,\dots,\mu_m;q)= \sn (w_1)\dots \sn(w_m) 
J_L(w_1(\mu_1),\dots,w_m(\mu_m);q).$$

It's more convenient to use the following normalization of
$J_L$:
$$Q_L(\mu_1,\dots,\mu_m;q)= 
J_L(\mu_1,\dots,\mu_m;q) J_U(\mu_1;q)\dots J_U(\mu_m;q),
$$
where $U$ is the unknot with  framing 0.
The invariant $Q_L$ enjoys the symmetry described in the following
proposition, which follows immediately from the definition.
\begin{pro} 

a) $Q_L(\mu_1,\dots,\mu_{n-1};q)$ is invariant
under the action of the Weyl group, i.e. for every $w_1,\dots,w_m\in\Weyl$, one has
$$Q_L(\mu_1,\dots,\mu_{m};q) = 
Q_L(w_1(\mu_1),\dots,w_{m}(\mu_{m});q).$$

b) If one of the $\mu_j$,\ $j=1\dots,m$, is on the boundary of the 
fundamental chamber $C$, then $Q_L(\mu_1,\dots,\mu_{m};q)=0$.
\end{pro}

For a positive integer $r$ let
$$F'_L(r)= \sum_{j=1}^m\sum_{\mu_j\in \Lambda'_r}Q_L(\mu_1,\dots, \mu_m;\zeta),$$
 where $\zeta^a = \exp(2a\pi i/r)$ for every rational number $a$, and
$$\Lambda'_r= \{\,\mu\in\Po \mid (\mu|\theta)\le r\,\} =\{\,
\sum_{i=1}^{n-1} k_i\lambda_i \in\Po \mid
 \sum k_i \le r\,\}.$$
Suppose the closed oriented 3-manifold $M^3$ is obtained from $S^3$ by surgery
 along a framed unoriented 
link $L$. Providing $L$ with arbitrary orientation,  the quantum $SU(n)$-invariant 
of $M$ at the $r$-th root of unity is defined by (see, for
example, \cite{KohnoTakata1,TuraevWenzl})
$$
\tau_r^{SU(n)}(M)= \frac{F'_L(r)}{F'_{U_+}(r)^{\sigma_+}\,
F'_{U_-}(r)^{\sigma_-}}.
$$

Here $\sigma_+,\sigma_-$ are the numbers of positive and negative
eigenvalues of the linking matrix of $L$ and  $U_\pm$ are the trivial knots
with framing $\pm 1$.

\subsection{Quantum invariants of links at roots of unity}
Let $\SSS$ be the simplex which is the convex hull of  the points 0, $r\lambda_1,\dots,r\lambda_{n-1}$. In other
words,
$$\SSS=\{\,z\in C \mid (z,\theta) \le r\,\} = \{ \sum_{i=1}^{n-1} a_i \lambda_i \in C\mid
 \sum_{i=1}^{n-1} a_i \le r\},$$
(we call it  the fundamental alcove of level $r$). Then $\Lambda'_r= C_r\cap \Po$. 
The  affine Weyl group at level $r$, by definition,  
is the group $\Weylr$ generated by all the reflections along boundary 
facets of the simplex $\SSS$.  It is known that  $\Weylr$ is the semi-direct
product of $\Weyl$ and the translation group $r\Lroot$ (see \cite{Kac}, chapter 6).

\begin{lem} \label{inv}
 Each of 
the lattices
$\Lroot$ and $\rho+\Lroot$ is invariant under the action of $\Weylr$,
for every positive integer $r$.
\end{lem}
\begin{pf} The root lattice $\Lroot$ is certainly invariant under the action of
the Weyl group and the translation group $r\Lroot$. Hence $\Lroot$ is
invariant under the action of $\Weylr$.

If $\alpha_i$ is a basis root, then the reflection along the
hyperplane perpendicular to $\alpha_i$ maps $\rho$ to $\rho-\alpha_i$,
which belongs to $\rho+\Lroot$. 
These reflections generate the Weyl group, hence  $\rho+\Lroot$ is invariant
under the action of the  Weyl group. It is obviously invariant
under the translation group $r\Lroot$.
\end{pf}

When specialized to the $r$-th root of unity, $Q_L$ has more symmetry.
The proof of the following important proposition is given in \S\ref{proof1}. 
\begin{pro} \label{91}

 Let $\zeta =\exp(2\pi i/r)$, where $r$ is a positive integer.

a)
$Q_L(\mu_1,\dots,\mu_{m};\zeta)$ is invariant
under the action of the affine Weyl group $\Weylr$, i.e. for every $w_1,\dots,w_m$
in   $\Weylr$, one has
$$Q_L(\mu_1,\dots,\mu_{m};\zeta) = 
Q_L(w_1(\mu_1),\dots,w_{m}(\mu_{m});\zeta).$$

b) If one of the $\mu_j, j=1\dots,m$, is on the boundary of the
 fundamental alcove  $\SSS$, then $Q_L(\mu_1,\dots,\mu_{m};\zeta)=0$. 
\end{pro}

\subsection{On  the definition of $PSU(n)$-quantum invariants}\label{defi}
Quantum
$PSU(n)$-invariants were introduced by Kirby-Melvin
\cite{KirbyMelvin} in the $n=2$ case, and  by
Kohno-Takata \cite{KohnoTakata2} in the $n>2$ case. 
We will give here another definition which is more natural from
our point of view. We will prove that our definition coincides with
that of Kohno and Takata in Appendix.

Let $\Lrootr$ be the set of all elements of the root lattice
 in the half-open parallelepiped
spanned by $r\alpha_1$, $\dots$, $r\alpha_{n-1}$, i.e. 
$$\Lrootr =\{\,k_1\al_1+\dots +k _{n-1}\al_{n-1}\in\Lroot 
 \mid 0\le k_i < r\,\}.$$
Define $F_L(r)$ using the same formula of $F'_L(r)$, only replacing 
$\Lambda'_r$ by $\rho+\Lrootr$:
$$F_L(r)= \sum_{j=1}^m\sum_{\mu_j\in(\rho+\Lrootr)}Q_L(\mu_1,
\dots, \mu_m;\zeta).$$

Suppose that $r$ and $n$ are coprime. We will show later that
$F_{U_\pm}(r)\not = 0$. Suppose that $M$ is obtained by surgery along the
framed link $L$.  Define the
quantum $PSU(n)$-invariant $\tauPSU_r(M)$ by

\begin{equation}\label{sao51}
\tauPSU_r(M)= \frac{F_L(r)}{F_{U_+}(r)^{\sigma_+}
F_{U_-}(r)^{\sigma_-}}.
\end{equation}

\begin{rem}
a) In the set $\rho+\Lrootr$ there are weights which are not
highest weights of any finite-dimensional $sl_n$-modules.

b) In the definition of $F_L$, the $\mu_j$'s run the set $\rho+\Lrootr$ since we 
use $V_{\mu}$ for the $sl_n$-module with highest weight $\mu-\rho$. Also if $n$ is
 {\em odd}, then $\rho$ is in the root lattice, hence $\rho+\Lroot=\Lroot$; and 
using the symmetry of $Q_L$ we can replace $\rho+\Lrootr$ by $\Lrootr$.

c) 
If in the definition of $F_L$, we let the $\mu_j$'s run the set
of all integral weights in the same parallelepiped, 
then from the above formula we get exactly the $SU(n)$-invariant. 
So both $SU(n)$ and $PSU(n)$-invariants can be defined by one formula, only for 
 the $SU(n)$ invariant, we use the weight lattice, while
for the  $PSU(n)$ invariant, we use the root lattice.
It is known that $\Lambda/\Lroot$ is the cyclic group of order $n$.

c) Kohno and Takata \cite{KohnoTakata2} showed that the $PSU(n)$ version is
finer than the $SU(n)$ version.
\end{rem}

The following important integrality property was established  by H. Murakami 
\cite{Murakami} for
$n=2$ , and by Takata-Yokota \cite{TakataYokota} and Masbaum-Wenzl
\cite{MasbaumWenzl} for $n>2$; 
this result can also be obtained by the method of this paper.

\begin{thm}
Suppose  $r$ is a prime not dividing $n!|H_1(M,\Z)|$.
Then the $PSU(n)$-quantum invariant
$\tauPSU_r(M)$ is a number in $\Z[\zeta]$.
Here $\zeta=\exp(2\pi i/r)$.

\end{thm}

\subsection{Existence of perturbative expansion}

Unlike the link case, quantum invariants of 3-manifolds can be
defined only at roots of unity.
In perturbative theory, we want to expand $q=e^h$, with $h$ an {\rm indeterminate}. 
However,
in the 3-manifold case, $q$ must be specified at the $r$-th root of unity, and 
it seems difficult to make meaning the substitution $q= e^h$, with
$h$ an indeterminate.
One way around is the following. Instead of $h$, we use $x=q-1$. Suppose $r$
is a prime not dividing $n!|H_1(M,\Z)|$.

Let $f(q)\in \Z[q]$ be a representative of $\tauPSU_r(M)\in \Z[\zeta]=
\Z[q]/(q^{r-1}+\dots+q+1)$.
Using the substitution $q=1+x$, we get

$$f(q)|_{q=x+1}= c_{r,0} + c_{r,1}x +\dots + c_{r,\ell} x^\ell+\dots.$$
The integer coefficients $c_{r,\ell}$ are, in general, dependent on
the representative $f(q)$. However, it's easy to prove the following (see \cite{Ohtsuki1}).
\begin{lem}
For $\ell{\le r-2}$, the classes $c_{r,\ell}\ \pmod{r}$ do
not depend on the representative $f(q)$
in $\Z[q]$.
\end{lem}
Hence for each prime
$r$, the classes
$c_{r,\ell} \pmod{r}$, $\ell=0,1,\dots,r-2$ are invariants of $M$. 
Let us fix $\ell$ and let $r$ vary (but $r$ must be an odd prime).
For each $r$ we have a residue class $c_{r,\ell}(M) \pmod r$.
We  want to show that 
these classes (when $r$ varies but $\ell$ fixed) can be unified in the following sense:
they are derived from  the same {\em rational
number} not depending on $r$. More precisely, (see Theorem \ref{main} below) 
we will show that there is a rational number $c_\ell(M)\in
\Z[\frac{1}{(2\ell+n(n-1))!\,|H_1(M,\Z)|}]$ such that 
\begin{equation}\label{444}
   c_\ell \equiv c_{r,\ell} \ \pmod{r}
\end{equation}
 for sufficiently
large prime $r$ (in fact, $r>2\ell+ n(n-1)$ is enough). 
The following is easy to prove (see \cite{Ohtsuki1}).
\begin{lem} \label{unique}
For each $\ell$, there exists at most  {\em one} rational number
$c_\ell\in \Z[\frac{1}{(2\ell+n(n-1))!\,|H_1(M,\Z)|}]$ such that
(\ref{444}) holds 
for all sufficiently large prime $r$.
\end{lem}

It follows that if such $c_\ell$ exist, then it is an invariant of
$M$. The existence of such $c_\ell$ is  a quite non-trivial fact.
If all the $c_\ell$'s do exist, then the series $\sum c_\ell x^\ell$  should be considered 
the perturbative expansion of $\tauPSU_r(M)$.

Another way to look at the existence of $c_\ell$ is the following.
For $f(q)\in\Z[q]$ (or in $\Zr[q]$) we will define ${\frak p}_r(f(q))$ as follows.
 First we substitute $q=x+1$ in $f$, then
cancel any powers of  $x$ with degree greater than 
$\bar r= (r-n(n-1)-1 )/2 $, finally reduce  all the coefficients 
modulo $r$. The result is ${\frak p}_r(f)$. 
So ${\frak p}_r$ is an algebra homomorphism
$${\frak p}_r \colon \Zr[q] \to \Z[x]/(r, x^{\bar r +1}).$$

It is easy to see that ${\frak p}_r$ descends to an algebra homomorphism, also denoted
by ${\frak p}_r$, on 
$$\Zr[q]/(q^{r-1}+\dots+q+1)=\Zr[\zeta].$$ 
Hence there is defined
${\frak p}_r(\tauPSU_r(M))$. Certainly
${\frak p}_r(x^d)=0$ if $d > \bar r$.

Recall that  $\Rin$ is the set of all
formal power series $\sum d_\ell x^\ell$ such that 
$d_\ell \in \Z[\frac{1}{(2\ell + n(n-1))!\,|H_1(M,\Z)|}]$. So, when $r$ is a prime
greater than both $|H_1(M,\Z)|$ and $2\ell + n(n-1)$, 
we can reduce $d_\ell$ modulo $r$.

For a series $g=\sum d_\ell x^\ell\in\Rin$, where
$|H_1(M,\Z)| < r$, 
we can also define ${\frak p}_r(g)$ in a similar way, as follows.
First cancel  powers of $x$ with degree greater than 
$\bar r$, then reduce  all the coefficients 
modulo $r$. The result is  ${\frak p}_r(g)$. 
 A similar operator  was also 
considered by Ohtsuki and Rozansky in connection with the $n=2$ case.

Then the equality $c_\ell \equiv c_{r,\ell}  \pmod{r}$, for every $\ell$ 
and $r> 2 \ell + n(n-1)$, 
means that 
$${\frak p}_r(\tauPSU_r(M))={\frak p}_r(\sum_\ell c_\ell x^\ell).$$
The following is the main result of this paper.

\begin{thm} \label{main} Suppose that $M$ is a rational homology 3-sphere.
There exists a  power series 
$$\tau^{PSU(n)}(M)= \sum_{\ell=0}^\infty c_\ell(M) x^\ell\,
 \in\, \Rin$$
such that for any prime number $r$ greater than
$\max \{n(n-1),|H_1(M,\Z)|\}$,
one has 

$${\frak p}_r(\tau^{PSU(n)}(M))={\frak p}_r(\tauPSU_r(M)).$$
Here $\tauPSU_r(M)$ is the $PSU(n)$-quantum invariant of $M$ at
the $r$-th root of unity.
\end{thm}

\begin{rem} 
a)   In the case $n=2$, the existence of $c_\ell$ was proved
by Ohtsuki \cite{Ohtsuki1}. Our normalization
of $\tauPSU$ differs from that of Ohtsuki by the factor $(\frac{|H_1(M,\Z)|}{r})$,
the Legendre symbol.

b) The above equality means that $\tau^{PSU(n)}(M)$ is
 the Fermat limit of $\tauPSU_r(M)$, in the terminology of \cite{LinWang}.
\end{rem}

The theorem shows that we can recover part of the quantum invariant
$\tauPSU_r(M)$ from the perturbative invariant $\tauPSU(M)$. 
We have the following conjecture.

\begin{conj}
The perturbative $PSU(n)$-invariant $\tau^{PSU(n)}(M)$ dominates
the quantum invariants $\tau_r^{SU(n)}(M)$, for every positive
integer $r$, not necessarily prime.
\end{conj}

Another conjecture is the integrality of $c_\ell(M)$, which is a generalization
of conjectures made and investigated by Lawrence, Lin, Rozansky and Wang
(see \cite{Lawrence,LawrenceRozansky,Rozansky2,LinWang}).
\begin{conj}\label{int}
The perturbative invariants $c_\ell(M)$ are  in $\Z[\frac{1}{n!\, |H_1(M,\Z)|}]$.
\end{conj}

One can show that $c_1(M) = (\frac{|H_1(M,\Z)|}{r})
 |H_1(M,\Z)|^{-n(n-1)/2}n(n^2-1) \lambda_C(M)$,
where $\lambda_C$ is the Casson-Walker invariant in Lescop normalization
\cite{Lescop}.

\section{Proof of the main theorem}
We will give the proof of the main theorem based on some technical
results which will be proved later.
\subsection{Integrality of the exponents} \label{183}
It is known that
for fixed $\mu_1,\dots,\mu_m$, the quantum invariant
$J_L(\mu_1,\dots,\mu_m;q)$ is a polynomial in $q^{1/2n}, q^{-1/2n}$ 
with {\em integer} coefficients (see, for example, \cite{MasbaumWenzl}
for a proof), i.e. $J_L(\mu_1,\dots,\mu_m;q)\in \Z[q^{1/2n},q^{-1/2n}]$. 
When the linking matrix of $L$ is 0, we have the following stronger result,
 whose proof will presented in \S\ref{proof1}.

\begin{pro}\label{31}
 Suppose that the linking matrix of
 $L$ is 0. Then $Q_L(\mu_1,\dots,\mu_m;q)$ is a polynomial
in $q, q^{-1}$ with integer coefficients. In other words, 
$Q_L(\mu_1,\dots,\mu_m;q)$ does not contain any non-integer
fractional power of $q$.
\end{pro}
\begin{rem} The proposition does not hold true if the linking matrix of
$L$ is not 0. 
\end{rem}
\begin{cor} Suppose that the linking matrix of $L'$ is {\em diagonal}
and that 
 $\mu_1,\dots,\mu_m$ are in
$(\rho+\Lroot)$. Then $Q_{L'}(\mu_1,\dots,\mu_m;q)$ is in 
$\Z[q,q^{-1}]$.\end{cor}

\begin{pf} Let $L$ be the same link as $L'$, only with 0 framing on each
component. By the previous proposition, $Q_L(\mu_1,\dots,\mu_m)$ is $\Z[q,q^{-1}]$.
From the general theory of
quantum invariants, it is known that  increasing  the framing by 1
on the $j$-th component results in a factor $q^{(|\mu_j|^2-
|\rho|^2)/2}$ in $Q_L$ (see, for example, \cite{Turaev}). 
Here $|\mu|^2= (\mu|\mu)$. We have:
\begin{equation}\label{1515}
Q_{L'}(\mu_1,\dots,\mu_m;q)= \prod_{j=1}^m q^{b_j\frac{|\mu_j|^2-|\rho|^2}{2}} 
Q_L(\mu_1,\dots,\mu_m;q).
\end{equation}
where $b_i$'s are the framing of components of $L'$.
It remains to notice that when $\mu$ is in $\rho+\Lroot$, $(|\mu|^2-|\rho|^2)/2$
is an integer.
 \end{pf}

The corollary shows an  advantage of using the root lattice  and links
with diagonal linking matrix (which correspond to rational homology 3-spheres).
If $\mu_j$ are in the bigger lattice $\Lambda$, then, in general, 
$Q_{L'}(\mu_1,\dots,\mu_m;q)$
contains fractional powers of $q$.

\subsection{Polynomial functions on the weight lattice}
Let us introduce a new  indeterminate $x=q-1$. Then 
$q^{-1}= 1-x + x^2-x^3+\dots$, 
and $Q_L$ becomes a formal power series in $x$:

\begin{equation}\label{sao21}
Q_L(\mu_1,\dots,\mu_m;q)\bigr|_{q=x+1}=\sum_{\ell=0}^\infty f_\ell(\mu_1,
\dots,\mu_m)x^\ell.
\end{equation}

This kind of expansion, actually,  a version of it for the $n=2$  case, was
first considered by Melvin and Morton \cite{MelvinMorton}.

Each monomial $\alpha_{i_1}\alpha_{i_2}\dots\alpha_{i_s}$, where
$\alpha_i$'s are the basis roots, defines a function on $\Lambda$ by
$$\alpha_{i_1}\alpha_{i_2}\dots\alpha_{i_s}(\mu)=(\alpha_{i_1}|\mu)
(\alpha_{i_2}|\mu)\dots(\alpha_{i_s}|\mu),$$
for $\mu\in\Lambda$.
Using linearity, every polynomial in the $\alpha_i$'s can be regarded as
a function on $\Lambda$; we call them the {\em polynomial functions} on $\Lambda$. 
The degree of the polynomial function is the degree of
the corresponding polynomial. A function of many variables, $f(\mu_1,\dots,\mu_m)$,
 with $\mu_1,\dots,\mu_m\in\Lambda$, is a
{\em polynomial function} if it is a linear combination of functions
of the form $g_1(\mu_1)\dots g_m(\mu_m)$, where each $g_j$ is a polynomial function. 

Let $\Delta$ be the following polynomial function on $\Lambda$ of degree $n(n-1)/2$:
$$\Delta =  \frac{\prod_{\alpha\in\Phi_+} \alpha} {\prod_{\alpha\in\Phi_+} (\rho|\alpha)}.$$
The denominator of 
$\Delta$ is $\prod_{\alpha\in\Phi_+} (\rho|\alpha)=
1^{n-1} 2^{n-2}\dots (n-1)^1$.
 
\begin{lem} \label{dimension}
a) If $\mu\in\Lambda_{++}$, then the dimension of $V_\mu$ is given by
$$\dim (V_\mu) = \Delta(\mu).$$

b) One has that $\Delta(\mu)=0$ if $\mu$ is on the boundary of
the the fundamental chamber $C$. For every $w\in\Weyl$,
$$ \Delta(w(\mu))= \sn(w) \Delta(\mu).$$
\end{lem}
Part a) is the famous dimension formula of Weyl. Part b) follows from elementary
Lie algebra theory (see \cite{Humphreys}).

The proof of the following proposition will be given in \S\ref{proof1}.
\begin{pro}\label{sao22}
Suppose that the link $L$ has $m$ components and $0$ linking matrix.

a) All the functions $f_i$ in (\ref{sao21}) are polynomial functions
divisible by $\Delta^2(\mu_1)\dots\Delta^2(\mu_m)$:
$$f_\ell(\mu_1,\dots,\mu_m)=\Delta^2(\mu_1)\dots\Delta^2(\mu_m) f'_\ell
(\mu_1,\dots,\mu_m).$$
Moreover, $f_\ell$ are polynomial functions taking {\em integer} values
when $\mu_1,\dots,\mu_m\in \Lambda$.

b)  Suppose the degree of $\mu_j$ in $f'_\ell$ is $s_j$. Then
$\sum_{j=1}^m s_j \le 3\ell/2$ and $\max_js_j \le 2\ell -\sum_{j=1}^m s_j$.
 For $\ell=0$
we have $f'_0=1$.
\end{pro}

\subsection{An expression of $Q_L$ in terms of exponential functions}
Let $q^\beta$, for $\beta\in\Lroot$, be the  formal power series in $x$
obtained by formally expanding  $q^\beta=(1+x)^\beta$,
the coefficients being polynomial in $\beta$:
$$q^\beta= \sum_{\ell=0}^\infty {\beta \choose \ell} x^\ell.$$ 
Here $\displaystyle{\beta\choose \ell}$ is the usual binome polynomial:
$\displaystyle{{\beta\choose \ell} = \frac{\beta(\beta-1)\dots(\beta-\ell+1)}{\ell!}}.$

Then $q^\beta(\mu)$, for $\beta\in\Lroot$ and
$\mu\in\Lambda$, is a formal power series in $x$ with integer coefficients:
$$q^\beta(\mu)= \sum_{\ell=0}^\infty {\beta \choose \ell}(\mu) x^\ell \in \Z[[x]],$$ 
since $(\beta|\mu)$ is an integer for every $\mu\in\Lambda$ and $\beta\in\Lroot$.

For each $i=1,\dots, n-1$, the element
$$\eta_i = 1-q^{-\al_i}$$
can be regarded as a  function from $\Lambda$ to $\Z[[x]]$. Moreover, the degree
0 term of $\eta_i(\mu)$ is 0, i.e. $\eta_i(\mu)$ is divisible by
$x$.

For $\ab=(a_1,\dots,a_{n-1})\in(\Zplus)^{n-1}$, let $|\ab|=a_1+\dots+a_{n-1}$, and
$\displaystyle{\boldeta^{\ab}= \prod_{i=1}^{n-1} \eta_i^{a_i}.}$
Then for $\mu\in\Lambda$, $\boldeta^\ab(\mu)\in x^{|\ab|} \Z[[x]]$. 

\begin{pro} Suppose that $L$ is a framed link  with $m$ components
and  0 linking matrix.  As formal power series in $x$ we have
\begin{equation}
\label{sao23}
Q_L(\mu_1,\dots,\mu_m;x+1)= \sum 
c_{\ab_1,\dots,\ab_m;\ell}\, \boldeta^{\ab_1}(\mu_1)\dots 
\boldeta^{\ab_m}(\mu_m) x^{\ell -\sum_j |\ab_j|},
\end{equation}
where the sum is over the set of $\ab_j\in (\Zplus)^{n-1}$ and $\ell\in \Zplus$ such that

\begin{equation}\label{cond0}
|\ab_j| \ge n(n-1)  \qquad \text {for every $j=1,\dots,m$},
\end{equation}
\begin{equation}\label{cond1}
\sum_{j=1}^m (|\ab_j|-n(n-1)) \le 3\ell/2, \quad \text{and}
\end{equation}
\begin{equation}\label{cond2}
\max_j (|\ab_j| - n(n-1)) \le 2\ell - \sum_j(|\ab_j| -n(n-1)).
\end{equation}
Moreover, the coefficient
$c_{\ab_1,\dots,\ab_m;\ell}$ is in $\Z_{(2\ell + (m+1) n(n-1)-\sum|\ab_j| +1)}$. 
 \label{987}
\end{pro}
\begin{rem} 
Since $\boldeta^{\ab}$ is a power series divisible by $x^{|\ab|}$,
the right hand side of (\ref{sao23})contains only non-negative powers of $x$.
The range of the sum is a bit complicated, the reader should just keep in mind
that $Q_L$ has the above form with appropriate restrictions on the indices
and on denominators of the coefficients.
\end{rem}
 
\begin{pf} The proof consists of a few changes of variables.
It is well known that if a polynomial $p(t_1,\dots,t_s)$
takes integer values whenever $t_i$'s are integers, then $p(t_1,\dots,t_s)$
is a $\Z$-linear combination of terms of the form
$\displaystyle{{t_1 \choose l_1} \dots {t_s\choose l_s}}$,
where $l_1,\dots,l_s$ are non-negative integers.

Recall that $(\alpha_i|\lambda_j)=\delta_{ij}$.
For $\ab = (a_1,\dots,a_{n-1})\in (\Zplus)^{n-1}$, let
$${{\boldsymbol\alpha} \choose \ab}= \prod_{i=1}^{n-1} {\alpha_i \choose a_i}.$$

By Proposition \ref{sao22}, $f_\ell(\mu_1,\dots,\mu_m)$ is a polynomial function
on $\mu_1,\dots,\mu_m$ which takes integer values whenever $\mu_j\in\Lambda$.
Hence 
\begin{equation}\label{s1s}
f_\ell(\mu_1,\dots,\mu_m) =\suml c'_{\ab_1,\dots,\ab_m;\ell}\, 
{{\boldsymbol\alpha}\choose {\ab_1}}(\mu_1)\dots 
{{\boldsymbol\alpha}\choose {\ab_m}}(\mu_m),
\end{equation}
with  $c'_{\ab_1,\dots,\ab_m;\ell}$ in $\Z$. 
Here $\suml$ means the sum is over the set of all $\ab_j\in(\Zplus)^{n-1}$
and $\ell \in \Zplus$ satisfying (\ref{cond0}), (\ref{cond1}) and (\ref{cond2})
These restrictions on $|\ab_j|$ and $\ell$ 
follows from the restriction on the degree of $f'_\ell$
(see Proposition \ref{sao22}, here $|\ab_j|= s_j + n(n-1)$).

Recall that $\eta_i= 1-(1+x)^{-\alpha_i}$. We can easily express $\alpha_i$ in terms
of  $x$ and $\eta_i$:
$$\alpha_i= \frac{\ln(1-\eta_i)}{\ln(1+x)},$$
and hence
\begin{equation}\label{919}
x \alpha_i= \frac{x \ln(1-\eta_i)}{\ln(1+x)}.
\end{equation}
Note that the right hand side can be expressed as a formal power series in
$x$ and $\eta_i$; moreoever, the coefficient of a term of total degree $d$
(in 
$x$ and $\eta_i$) 
is in $\Z_{(d+1)}$.

Using (\ref{s1s}) in formula (\ref{sao21}), we get
\begin{align}
Q_L(\mu_1, \dots,\mu_m;q)\bigr|_{q=x+1}& =\sum_{\ell=0}^\infty 
\suml 
c'_{\ab_1,\dots,\ab_m;\ell}\, x^\ell \prod_{j=1}^m
{{\boldsymbol\alpha}\choose {\ab_j}}(\mu_j)
 \notag \\
&= \sum_{\ell=0}^\infty 
\suml 
c'_{\ab_1,\dots,\ab_m;\ell}\, x^{\ell-\sum_j|\ab_j|}\prod_{j=1}^m 
{{\boldsymbol\alpha}\choose {\ab_j}}(\mu_j) \, x^{|\ab_j|}.
\notag 
\end{align}
Using the change of variable (\ref{919}) in the above formula, with
some simple calculation,  we get (\ref{sao23}).
The denominators of the coefficients come from the denominators
of the expressions ${{\boldsymbol \alpha} \choose \ab_j}$ and
from the denominators of the right hand side of (\ref{919}).
It is easy to check that $c_{\ab_1,\dots,\ab_m;\ell}$ is
in $\Z_{(2\ell + (m+1) n(n-1)-\sum|\ab_j| +1)}$.
\end{pf}

\subsection{Linear operator $\Gamma_b$} We now introduce the linear operator
$\Gamma_b$ which is a kind of Laplace transform.
For $\mu\in\Lambda$ let
\begin{equation}\label{psi}
\psi(\mu;q)= \prod_{\alpha\in\Phi_+}(1-q^{-(\mu|\alpha)})
=\prod_{\alpha\in\Phi_+}(1-q^{-\alpha})(\mu) .
\end{equation}
This function plays important role in Lie theory, it appeared in  Weyl's
character formula.

Let $\Sa$ be the free $\Z$-module generated by $q^\beta$, with $\beta\in\Lroot$.
For each  integer $b$ not divisible by the odd prime $r$ let $\Gamma_b$ be the linear operator
$$\Gamma_b \colon \Sa \to \Z[\frac{1}{n!b}][[x]], \qquad \text{defined by}$$
$$\Gamma_b(q^{\beta})=  (1+x)^{(-|\beta|^2/2b)}\, y_b.$$
Here $y_b$ is an element in $\Z[\frac{1}{n!}][[x]]$, not depending
on $\beta$, and is given by
$$y_b = \frac{1}{n!}\left(\frac{|b|}{r}\right)^{n-1}
(1+x)^{\frac{\sn(b)-b}{2} |\rho|^2}\, \psi\left(-\sn(b)\rho;1+x
\right),$$
($sn(b)$ is the sign of $b$).

Note that $\boldeta^{\ab}$ is in $\Sa$, hence there is defined $\Gamma_b(\boldeta^{\ab})$.
A property of $\Gamma_b$ is stated in the following proposition, whose
proof will be presented in \S\ref{gammb}.

\begin{pro} \label{new1}
 The power series $\Gamma_b(\boldeta^\ab)$ is divisible by
$x^{\lfloor (n(n-1) + |\ab|+1)/2\rfloor}$. Here $\lfloor z\rfloor$ is the
integer part of $z$.
\end{pro}

\subsection{Formula of perturbative invariants $\tauPSU$}
We continue to assume that $L$ has 0 linking matrix with $m$ components.
Suppose that $L'$ is the same as $L$, except
that the framing on components are {\em non-zero} integers $b_1,\dots,b_m$. Suppose for
$L$ we have the expansion (\ref{sao23}) of $Q_L$.

Let $\tau(L')$ be obtained from $Q_L$ (i.e. the right hand side of (\ref{sao23}))
by replacing $\boldeta^{\ab_j}(\mu_j)$ by
$\Gamma_{b_j}(\boldeta^{\ab_j})$:

\begin{equation}\label{988}
\taut(L')= \sum c_{\ab_1,\dots,\ab_m;\ell}\, 
\Gamma_{b_1}(\boldeta^{\ab_1})\dots\Gamma_{b_m}(\boldeta^{\ab_m})\,
 x^{\ell -\sum_j |\ab_j|}.
\end{equation}
The range of the sum is the same as in formula (\ref{sao23}), i.e.
the sum is over the set of all $\ab_j\in\Zplus^{n-1}$ and $\ell\in\Zplus$
satisfying (\ref{cond0}), (\ref{cond1}) and (\ref{cond2}).

By  Proposition \ref{new1}, the term 
$\Gamma_{b_1}(\boldeta^{\ab_1})\dots\Gamma_{b_m}(\boldeta^{\ab_m})
 x^{\ell -\sum_j|\ab_j|}$ is divisible by the power of $x$ with exponent 
$$\ell -\sum_j|\ab_j| + \frac{mn(n-1)}{2} + \sum_j 
\lfloor \frac{|\ab_j|}{2}\rfloor,$$
which is greater than or equal to $\ell -\sum_j |\ab_j|/2 + mn(n-1)/2$.
The latter, by (\ref{cond2}), is greater than or equal to
$\frac{1}{2}\max_j(|\ab_j|-n(n-1))$. For each fixed number 
$a$, there are only a finite number of $\ab_j\in(\Zplus)^{n-1}$ such that
$\frac{1}{2}\max_j(|\ab_j|-n(n-1)) < a$. 
Hence the right hand side of
(\ref{988}) is really a formal power series in $x$.

The denominators of coefficients of $\tau(L')$ come from the factor
$n!b_j$ and the denominators of $c_{\ab_1,\dots,\ab_m;\ell}$. Hence, using
the property of the denominators of $c_{\ab_1,\dots,\ab_m;\ell}$ in 
Proposition \ref{987}
 one can  easily shows
that in
$\tau(L')$ the coefficient of $x^\ell$ is in $\Z[\frac{1}{(2\ell+n(n-1))! b_1\dots
b_m} ]$. In other words, $\tau(L')$ is in $\Q[[x]]_M$, where 
$M$ is the 3-manifold obtained by surgery along $L'$.

We will  show that $\taut(L')$ is an invariant of $M$, which is
a rational homology 3-sphere.
Later in \S\ref{len} we will  show  that the constant term of 
the power series $\taut(L')$ is 
$$\left(\frac{|H_1(M,\Z)|}{r}\right)|H_1(M,\Z)|^{-n(n-1)/2}.$$ 
Hence there always exists the inverse $\taut(L')^{-1}\in\Rin$.

Let $M$ be an arbitrary  rational homology 3-sphere. Ohtsuki
showed (see \cite{Ohtsuki1}, see also \cite{Murakami})
that there are lens spaces $M(d_1),\dots, M(d_s)$ such that
$M\# M(d_1) \#\dots \# M(d_s)$ can be obtained from $S^3$ by
surgery on a link $L'$ with diagonal linking matrix. 
Here $M(d)$ is the lens space obtained by surgery along $U_d$, 
the unknot with framing $d$. Moreover each $d_i$ is less than or equal to
$|H_1(M;\Z)|$.
Define
$$\tau^{PSU(n)}(M)= \taut (L')\, \taut(U_{d_1})^{-1}\dots 
\taut(U_{d_s})^{-1}.$$

We have not proved that $\tau^{PSU(n)}(M)$ is an invariant of $M$ yet.
So, for the time being, $\tau^{PSU(n)}(M)$ means any formal power series
obtained by the above procedure, which a priori depends on the choice of
$d_1,\dots,d_s$ and $L'$.

The main theorem can be formulated in a more precise form as follows.

\begin{thm} \label{main1}
Suppose $M$ is rational homology 3-sphere, and $r$ a prime number
greater than $\max\{n(n-1),|H(M,\Z)|\}$. Then 
$${\frak p}_r\left(\tauPSU_r(M)\right)= {\frak p}_r\left(\tau^{PSU(n)}(M)\right).$$
\end{thm}
Now the uniqueness of Lemma \ref{unique}, with this theorem, shows that
 $\tau^{PSU(n)}(M)$ is really
an invariant of the 3-manifold $M$, i.e. it does not depend on the
 choice of $d_1,\dots,d_s$ and $L'$ in the construction of $\tau^{PSU(n)}(M)$.

Since the series $\Gamma_b(\boldeta^\ab)$ has coefficients in $\Z[\frac{1}{n!b}]$, we
see that if all the coefficients $c_{\ab_1,\dots,\ab_m;\ell}$ are in 
$\Z[\frac{1}{n!}]$,
then the power series $\tauPSU (M)$ has coefficients in 
$\Z[\frac{1}{n!|H_1(M,\Z)|}]$. Hence Conjecture \ref{int} follows from the following
conjecture.
\begin{conj}
All the coefficients $c_{\ab_1,\dots,\ab_m;\ell}$ are in $\Z[\frac{1}{n!}]$.
\end{conj}
 This conjecture, for the special case $n=2$ and in a slightly different form,
 had been formulated by Rozansky.
He also proved that the conjecture holds true in the special case $n=2$ and
the link $L$ has {\em one} component (a knot), using explicit formula
of the $R$-matrix (see \cite{Rozansky2}).

\subsection{Proof of the main theorem} 
Case 1: 
$M$ is obtained by surgery along
a framed link $L'$ with diagonal linking matrix and $\tau^{PSU(n)}(M)=\taut(L')$.
Let $L$ be the same link as $L'$  except that the framing of each
component is 0. Suppose that the
framing on components  of $L'$ are $b_1,\dots,b_m$. We will 
fix an odd prime number $r$ greater than $\max\{n(n-1), |H_1(M,\Z)|\}$. 
Let us consider the expression (\ref{sao23})
of $Q_L(\mu_1,\dots,\mu_m;x+1)$, 
and let $\Qone$ be the part of $Q_L$ in (\ref{sao23}), where the indices
satisfy 
\begin{equation} \label{cond3}
2\ell - \sum_j(|\ab_j|-n(n-1)) \le r - n(n-1).
\end{equation}
In other words,
\begin{equation}\label{one1}
\Qone(\mu_1,\dots,\mu_m;x+1)= \sum{}^{(1)}
c_{\ab_1,\dots,\ab_m;\ell}\, \boldeta^{\ab_1}(\mu_1)\dots 
\boldeta^{\ab_m}(\mu_m) x^{\ell -\sum_j |\ab_j|},
\end{equation}
where $\sum{}^{(1)}$ means the sum over the set of all
$\ab_j\in\Zplus^{n-1}$ and $\ell\in\Zplus$ satisfying 
(\ref{cond0}), (\ref{cond1}), (\ref{cond2})
and (\ref{cond3}). 
The reason we choose the new restriction (\ref{cond3}) is that, 
by Proposition \ref{987}, all the coefficients $ c_{\ab_1,\dots,\ab_m;\ell}$ 
in (\ref{one1}) are in
$\Zr$. Note that $\sum{}^{(1)}$ is a finite sum.

\begin{lem} \label{1414}
There is the following splitting for $Q_L$ (with $q= x+1$):
$$Q_L(\mu_1,\dots,\mu_m;q)= \Qone(\mu_1,\dots,\mu_m;q) +
\Qtwo(\mu_1,\dots,\mu_m;q)
+ \Qthree(\mu_1,\dots,\mu_m;q),$$
in which :

a) $\Qone$ is defined by (\ref{one1}),

 b) $\Qtwo$ is a polynomial in $(q-1)=x$ of the form
\begin{equation}\label{sao61}
 \Qtwo(\mu_1,\dots,\mu_m;q)= \sum{}^{(2)}
\tilde c_{\ab_1,\dots,\ab_m;\ell}{{\boldsymbol\alpha}\choose \ab_1}(\mu_1)
\dots {{\boldsymbol\alpha}\choose \ab_m}(\mu_m)\,
x^\ell,
\end{equation}
where $\sum{}^{(2)}$ is the sum over the set of all $\ab_j\in\Zplus^{n-1}$
and $\ell\in\Zplus$ satifying (\ref{cond0}),
(\ref{cond1}), (\ref{cond2}), and , in addition,
\begin{equation}\label{cond4}
 2\ell - \sum_j(|\ab_j|-n(n-1)) > r - n(n-1) ,
\end{equation}
\begin{equation}\label{cond5}
\ell < N.
\end{equation}

c) 
$\Qthree(\mu_1,\dots,\mu_m;q)$ is a polynomial in $q,q^{-1}$
which is divisible by $(q-1)^{N}$ in $\Zr[q,q^{-1}]$.

Here $N= m (r-n-1)(n-1)/2 + \bar r +1$.
\end{lem}
Note that (\ref{cond3}) and (\ref{cond4}) are complementary.
\begin{pf}
With the substitution $q=x+1$ and
\begin{equation}\label{sao24}
\eta_i(\mu)= 1 -(1+x)^{-(\alpha_i|\mu)},
\end{equation}
 $\Qone$ becomes a  formal power series in $x$.
The right hand side of (\ref{sao24}) is a power series with {\em integer}
coefficients, hence the denominators of the power series $\Qone$ comes from
the coefficients $c_{\ab_1,\dots,\ab_m;\ell}$ which are in $\Zr$. Hence
$\Qone$ is a power series in $x$ with coefficients in $\Zr$.
On the other hand, $Q_L$ is a power series in $x$ with integer coefficients.
Hence, for fixed $\mu_1,\dots,\mu_m$,
 $Q_L-\Qone$ is a power series in $x$ with coefficients in $\Zr$.

It follows that 
\begin{equation}\label{a1a}
Q_L-\Qone = \sum
\tilde c_{\ab_1,\dots,\ab_m;\ell}{{\boldsymbol\alpha}\choose \ab_1}(\mu_1)
\dots, {{\boldsymbol\alpha}\choose \ab_m}(\mu_m)
x^\ell,
\end{equation}
where the sum is over the set of all $\ab_j\in\Zplus^{n-1}$ and $\ell\in\Zplus$
satisfying (\ref{cond0}), (\ref{cond1}), (\ref{cond2}) and (\ref{cond4}), and
the coefficients $\tilde c_{\ab_1,\dots,\ab_m;\ell}$ are in
$\Zr$.

Let us define $\Qtwo$ by using the right hand side of the (\ref{a1a}),
imposing additional restriction  (\ref{cond5}) on the indices $\ab_j,\ell$.
It is easy to see that there are only a finite number of terms in
$\Qtwo$. Finally, let $\Qthree = Q_L -\Qone -\Qtwo$. A term of the
right hand side of (\ref{a1a}), for which (\ref{cond5}) is
not satisfied, must be divisible by
the power of $x$ with exponent $\ell+1$ which is
greater than or equal to $N$. Hence $\Qthree$ is divisible by
$x^N$.

Let us return to the indeterminate $q=x+1$. The advantage of
using this indeterminate is that all 
$Q_L,\Qone,\Qtwo,\Qthree$ are {\em polynomials} in $q,q^{-1}$
(i.e. they don't have infinite order), which can be seen as follows.

First of all,  $Q_L$ is a polynomial
in $q$ and $q^{-1}$ with coefficients in $\Z$, by Proposition \ref{31}. 

For fixed $\mu_1,\dots,\mu_m$, replacing $x=q-1$ in the expression
of $\Qtwo$, we see that $\Qtwo$ is a polynomial in $q$ with
coefficients in $\Zr$.

Note that, for $i=1,\dots, n-1$, the number $(\alpha_i|\mu)$ is always
an integer for every $\mu\in\Lambda$. Hence 
 $\eta_i(\mu)=1-q^{-(\alpha_i|\mu)}$ is a polynomial 
in $q,q^{-1}$ with integer coefficients,  and $\eta_i$
is divisible by $(q-1)$ in $\Z[q, q^{-1}]$. Now formula (\ref{one1})
shows that 
$\Qone$ is a polynomial in $q,q^{-1}$ with coefficients in $\Zr$.

It follows that the remaining  
$\Qthree$ is also a polynomial in $q,q^{-1}$ with coefficients in
$\Zr$. Moreover, $\Qthree$ is divisible by $(q-1)^{N}$ in
$\Zr[q,q^{-1}]$:

$$\Qthree(\mu_1,\dots,\mu_m;q)= (q-1)^{N}\,\sum_{s=-d}^{d} g(\mu_1,\dots,\mu_m) q^s,$$
for some positive integer $d$. Here $g(\mu_1,\dots,\mu_m)$, possibly not  a
polynomial function, takes values in $\Zr$
when $\mu_1,\dots,\mu_m\in\Lambda$. This completes the proof of the proposition.
\end{pf}

Recall that $L'$ differs from $L$ only by the framing on 
components. 
By formula (\ref{1515}):
$$
{L'}(\mu_1,\dots,\mu_m;q)= \prod_{j=1}^m q^{b_j\frac{|\mu_j|^2-|\rho|^2}{2}} 
Q_L(\mu_1,\dots,\mu_m;q),
$$
and hence
\begin{equation}\label{sao81}
F_{L'}(r)= \sum_{\mu_j\in\ (\rho+\Lrootr)}
\prod_{j=1}^m q^{b_j\frac{|\mu_j|^2-|\rho|^2}{2}} Q_L(\mu_1,\dots,\mu_m;\zeta).
\end{equation}

Let $\Fone(r),\Ftwo(r),\Fthree(r)$ be the value of the right hand side of (\ref{sao81}), if
we replace $Q_L$ by, respectively, $\Qone,\Qtwo,\Qthree$. Then we have
$$F_{L'}=\Fone +\Ftwo +\Fthree, \qquad \text{and}$$
\begin{align}
\tauPSU_r(M)&= \frac{F_{L'}(r)}{F_{U_+}(r)^{\sigma_+} F_{U_-}(r)^{\sigma_-}} \notag\\
&= \tauone + \tautwo +\tauthree,\notag \end{align}
where
$$\tau^{(j)}= 
\frac{F^{(j)}_{L'}(r)}{F_{U_+}(r)^{\sigma_+} F_{U_-}(r)^{\sigma_-}}.$$
The following 3 lemmas prove the main theorem (Theorem \ref{main1}) in case 1.
\begin{lem}\label{three} The number 
$\tauthree$ is divisible by $(\zeta-1)^{\bar r +1}$ in
$\Z_{(r)}[\zeta]$. Hence 
${\frak p}_r(\tauthree)=0.$
\end{lem}
\begin{lem}\label{two}
The number 
$\tautwo$ is divisible by $(\zeta-1)^{\bar r +1}$ in
$\Z_{(r)}[\zeta]$. Hence 
${\frak p}_r(\tautwo)=0.$

\end{lem}
\begin{lem}\label{one}
One has that \qquad
${\frak p}_r (\tauone)={\frak p}_r(\taut(L')).$
\end{lem}

The expression $F_{U_\pm}(r)$ appears in the denominators
of the formulas of
$\tau^{(j)}$. We will need the following proposition whose proof
will be given in \S\ref{proof2}.
\begin{pro} \label{sao82}
Each of $F_{U_+}(r)$ and $F_{U_-}(r)$ is
proportional to $(\zeta-1)^{(n-1)(r-n-1)/2}$ by a proportional factor
which is a unit in
$\Z_{(1)}[\zeta] = \Z[\frac{1}{n!}][\zeta]$.
\end{pro}

\begin{pf} [Proof of Lemma \ref{three}]
Since each $b_j$ is non-zero, $\sigma_+ +\sigma_-= m$.
By Proposition \ref{sao82}, $F_{U_+} (r)^{\sigma_+} F_{U_-}(r)^{\sigma_-}$
is proportional to $(\zeta-1)^{ m(n-1)(r-n-1)/2}$ by a unit 
in $\Zr[\zeta]$.
But $\Fthree(r)$ is divisible
by $(\zeta-1)^N$ in $\Zr[\zeta]$ by Lemma \ref{1414}. Hence
$$\tauthree= \frac{\Fthree(r)}{F_{U_+} (r)^{\sigma_+} F_{U_-}(r)^{\sigma_-}}$$
 is divisible by the power of $\zeta-1$ with exponent
$N-m(n-1)(r-n-1)/2=
\bar r +1.$
\end{pf}

To prove Lemma \ref{two} we will need the following lemma whose proof
will be given in \S\ref{proof2}.
\begin{lem}\label{sao71} Let $\ab\in(\Zplus)^{n-1}$ and $b\in\Zplus$.
Then the number
$$\sum_{\mu\in\rho +\Lrootr}
\zeta^{b\frac{|\mu|^2-|\rho|^2}{2}} {{\boldsymbol \alpha}\choose \ab}(\mu)$$
is proportional to $(\zeta-1)^{(n-1)(r-1)/2 - \lfloor |\ab|/2 \rfloor}$ 
by a factor in $\Zr[\zeta]$.
\end{lem}
\begin{pf}[Proof of Lemma \ref{two}]
Using the expression (\ref{sao61})
we have that
\begin{align}
\Ftwo(r)&= 
\sum_{\mu_j\in\rho +\Lrootr} \prod_{j=1}^m 
\zeta^{b_j\frac{|\mu_j|^2-|
\rho|^2}{2}}
\sum{}^{(2)} \tilde c_{\ab_1,\dots,\ab_m;\ell}
{{\boldsymbol \alpha}\choose \ab_1}(\mu_1)\dots
{{\boldsymbol \alpha}\choose \ab_m}(\mu_m) \, x^\ell
\notag \\
&=
\sum{}^{(2)} \tilde c_{\ab_1,\dots,\ab_m;\ell}\,x^\ell
\prod_{j=1}^{m} \sum_{\mu_j\in\rho +\Lrootr}
\zeta^{b_j\frac{|\mu_j|^2-|
\rho|^2}{2}} {{\boldsymbol \alpha}\choose \ab_j}(\mu_j)
\notag
\end{align}
 Applying Lemma \ref{sao71}, we see that
$\Ftwo(r)$ is divisible, in $\Zr[\zeta]$, by 
the power of $\zeta-1$ with the exponent
$$ \ell + m\frac{(n-1)(r-1)}{2} - \sum_j\lfloor\frac{|\ab_j|}{2}\rfloor,$$
which,  by (\ref{cond4}), is greater than or equal to $N$.
As in the proof of Lemma \ref{three}, it follows that
$\tautwo$ is divisible by $(\zeta-1)^{\bar r +1}$.
\end{pf}

In order to prove Lemma \ref{one}, we will need the following proposition,
whose proof will be presented in \S\ref{proof2}.
\begin{pro}\label{new2}
Let
$$g_b(\ab) = \frac
{\displaystyle{\summu}\zeta^{b\frac{|\mu|^2-|\rho|^2}{2}}
\boldeta^{\ab}(\mu)}
{F_{U_{\sn(b)}}(r)}.$$

Then \qquad ${\frak p}_r (g_b(\ab))= {\frak p}_r(\Gamma_b(\boldeta^{\ab}))$.
\end{pro}

\begin{pf}[Proof of Lemma \ref{one}].
Recall that 
$$\Qone = \sum{}^{(1)}
 c_{\ab_1,\dots,\ab_m;\ell}\,\boldeta^{\ab_1}(\mu_1)
\dots\boldeta^{\ab_m}(\mu_m)\,
(q-1)^{\ell - \sum_j|\ab_j|}.$$
Here $c_{\ab_1,\dots,\ab_m;\ell}\in \Zr$.
Hence

$$\Fone(r)= \sum{}^{(1)}
 c_{\ab_1,\dots,\ab_m;\ell} (\zeta-1)^{\ell-\sum_j|\ab_j|}
\prod_{j=1}^m
\sum_{\mu_j\in(\rho+\Lrootr)}  \zeta^{b_j\frac{|\mu_j|^2-
|\rho|^2  }{2}}\boldeta^{\ab_j}(\mu_j),$$
and
\begin{align}
\tauone & = \frac{\Fone(r)}{F_{U_+}(r)^{\sigma_+} F_{U_-}(r)^{\sigma_-}}
\notag\\
&=  \sum{}^{(1)} c_{\ab_1,\dots,\ab_m;\ell}
(\zeta-1)^{\ell -\sum_j|\ab_j|}
\prod_{j=1}^m
g_{b_j}(\ab_j). \notag
\end{align}
By Proposition \ref{new2},
${\frak p}_r (g_b(\ab))= {\frak p}_r(\Gamma_b(\boldeta^{\ab}))$. It follows that
$${\frak p}_r(\tauone)= {\frak p}_r\left( \sum{}^{(1)}
 c_{\ab_1,\dots,\ab_m;\ell}\, x^{\ell -\sum_j|\ab_j|}
\prod_{j=1}^m \Gamma_{b_j}(\boldeta^{\ab_j}) \right).$$
The argument of ${\frak p}_r$ on the right hand side of this formula, by
definition, is a part of $\taut(L')$ (see formula (\ref{988}) for
the definition of $\taut(L'))$. This part differs from
 $\taut(L')$  only by elements of degree $> \bar r$, which
are annihilated by ${\frak p}_r$, hence
$${\frak p}_r(\tauone)={\frak p}_r(\taut(L')).$$
This completes the proof of Lemma \ref{one}, and hence that of Main theorem in case 1.
\end{pf}

Case 2: $M$ is an arbitrary rational homology 3-sphere, 
 $M'=
M\#M(d_1)\#\dots \# M(d_s)$ is obtained by surgery along a link
$L'$ with diagonal linking matrix, and
$$\tau^{PSU(n)}(M) = \taut(L') \taut(U_{d_1})^{-1} \dots
\taut(U_{d_s})^{-1}.$$ 
The multiplicative property of quantum invariants says that

\begin{equation}\label{sao83}
\tauPSU_r(M')= \tauPSU_r(M)\, \tauPSU_r(M(d_1))\dots\tauPSU_r(M(d_s)).
\end{equation}
Applying the result of case 1 to $L'$ and $U_{d_1},\dots,U_{d_s}$, and
using (\ref{sao83}) we see that 
${\frak p}_r(\tauPSU_r(M))={\frak p}_r(\tau^{PSU(n)}(M))$. 
This completes the proof of the main theorem.

\section{Quantum invariants of links and the Kontsevich integral}
\label{proof1}
In this section we always assume that $L$ is a framed oriented link
of $m$ components, and that $r$ is a positive integer. We will investigate the 
dependence of the quantum invariant $J_L(\mu_1,\dots,\mu_m;q)$ on the $\mu_j$'s 
and  prove Propositions
\ref{91}, \ref{31}, and \ref{sao22}. 

\subsection{Some general facts about quantum invariants}
Recall that $J_L(\mu_1,\dots,\mu_m;q)$, or $J_L(V_{\mu_1},\dots,
V_{\mu_m};q)$ is the quantum invariant of $L$, when the components of
$L$ are colored by the $sl_n$-modules $V_{\mu_1},\dots,V_{\mu_m}$. Here
$V_\mu$ is the irreducible module of highest weight $\mu-\rho$, for
$\mu\in\Poo$. 

Recall that
for  $\mu\in\Lambda$, 
$$\psi(\mu;q)= \prod_{\alpha\in\Phi_+}(1-q^{-(\mu|\alpha)})
=\prod_{\alpha\in\Phi_+}(1-q^{-\alpha})(\mu) .$$
From the Weyl denominator formula (see, for example, \cite{Kac}),
 we have the following.

\begin{lem}\label{weyl}
One has that \qquad $\displaystyle
{\psi(\mu;q) = \sum_{w\in\Weyl} \sn(w) q^{(\mu|w(\rho)-\rho)}}$. 
\end{lem}

The value of the quantum invariant $J_L$ of the unknot $U$ is 
 (see, for example
 \cite{Jones}):
\begin{equation}\label{unknot}
J_U(\mu;q)= q^{(\mu-\rho|\rho)}\frac{\psi(\mu;q)}{\psi(\rho;q)} .
\end{equation}
Note that $(\mu|\alpha)\in\Z$ for $\mu\in\Lambda$ and $\alpha\in\Lroot$.
 Moreover, it's known that $2\rho\in\Lroot$, and hence
$(\mu|\rho)\in \frac{1}{2}\Z$. It follows that $J_U(\mu;q)$ is either
in $ q^{1/2}\Z[q,q^{-1}]$ or in $\Z[q,q^{-1}]$.

Let $L^{(2)}$ be the framed link obtained from $L$ by replacing  the first
component by two its parallel push-offs (using the frame). Suppose that
the color of these 2 push-offs are modules $V$ and $V'$. Then we have 
the following formula (see \cite{Turaev}):
$$
J_{L^{(2)}}(V,V',\dots;q)= J_L(V\otimes V',\dots;q).
$$
Hence for $Q_L$ we have a similar result:
\begin{equation}\label{parallel}
Q_{L^{(2)}}(V,V',\dots;q)= Q_L(V\otimes V',\dots;q).
\end{equation}

In particular, when $L=U$, noting that $U^{(2)}$ is the trivial link with
2 components, we have that
\begin{align}
 J_U(V\otimes V';q)&= J_{U^{(2)}}(V,V';q) \notag\\
&= J_U(V;q) J_U(V';q)\notag
\end{align}
This means
$$J_U(\cdot;q): \RR \to \Z[q^{\pm 1/2n}]$$
is a {\em ring} homomorphism. Recall that $\RR$ is the ring of finite-dimensional
$sl_n$-modules.
Let $I_r$ be the ideal in $\RR$ generated by all the $V_\mu$, with $(\mu|\theta)=r$, 
(here $\mu\in\Poo$). Recall that $\theta=\alpha_1+\dots+\alpha_{n-1}$ is the
longest root.
The quotient ring $\RR/I_r$ is known as the fusion algebra (see
\cite{GoodmanWenzl,KohnoTakata1} and reference therein).

\begin{lem} If $\mu$ is in $I_r$, then $J_U(\mu;\zeta)=0$. 
\end{lem}
\begin{pf} Since $J_U(\cdot;\zeta): \RR \to \C$ is a ring homomorphism, 
it suffices to show that
$J_U(\mu;\zeta)=0$ if $(\mu|\theta)=r$. 
We have that
$$1 -\zeta^{-(\mu|\theta)}= 1-\zeta^{-r}=0.$$
Since $\theta$ is one of the positive roots, it follows from the definition 
of $\psi$ that $\psi(\mu;\zeta)=0$, and hence
$J_U(\mu;\zeta)=0$.
\end{pf}

\begin{pro} Suppose that one of $\mu_1,\dots,\mu_m$ is in $I_r$.
Then $J_L(\mu_1,\dots,\mu_m;\zeta)=0$.
\label{101}
\end{pro}
\begin{pf}  A proof for the case $n=2$ is
given in \cite{KirbyMelvin}, Lemma 3.29. 
Suppose $\mu_1\in I_r$. The argument of Lemma 3.29 of
 \cite{KirbyMelvin} shows that $J_L(\mu_1,\dots,\mu_m;\zeta)$ is
proportional to $J_U(\mu_1;\zeta)$ which is 0. Hence 
$J_L(\mu_1,\dots,\mu_m;\zeta)=0$.
\end{pf}
\begin{rem} This proposition also follows from the theory of the Kontsevich integral,
more precisely, from Proposition \ref{Kont} and formula (\ref{2222}) below.
\end{rem}

\subsection{Proof of Proposition \ref{91}}
 We will fix $\mu_2,\dots,\mu_m$ and will
write $J_L(\mu_1;q)$ instead of $J_L(\mu_1,\mu_2,\dots,\mu_m;q)$.

Proposition \ref{101} shows that  $J_L(\mu_1;\zeta)=J_L(\mu'_1;\zeta)$ if
$\mu_1 \equiv \mu'_1  \pmod{I_r}$.

From the theory of fusion algebra it
is known that for every $w\in\Weylr$, 
$$ V_\mu \:\equiv\:  \sn (w) V_{w(\mu)}  \pmod{I_r}.$$
Hence 
$$J_L(\mu;\zeta)= \sn(w)J_L(w(\mu);\zeta).$$
Put $L=U$ we get 
$$J_U(\mu;\zeta)= \sn(w)J_U(w(\mu);\zeta).$$
Take the product of the last two identities, we get
$$Q_L(\mu;\zeta)= Q_L(w(\mu);\zeta),$$
which proves part a) of Proposition \ref{91}.
 Part b) follows from a) and Proposition
\ref{101}.

\subsection{Proof of proposition \ref{31}}

We have to show that $Q_L(\mu_1,\dots,\mu_m;q)\in\Z[q^{\pm 1}]$ if 
the linking matrix of $L$ is 0.  It's enough to consider the
case when $\mu_j\in\Poo$, since  $Q_L$ is invariant
under the action of the Weyl group and equal to 0 if one of the $\mu_j$ is
on the boundary of the fundamental chamber $C$.

The invariant $Q_L$ satisfies
the doubling formula (\ref{parallel}), and the representation
ring $\RR$ is generated (as algebra over $\Z$) by 
$V_{\lambda_1+\rho},\dots,V_{\lambda_{n-1}+\rho}$. Hence it suffices
to prove that $Q_L(\mu_1,\dots,\mu_m;q)\in\Z[q^{\pm 1}]$
 for the case when each $\mu_j$ is in
$\{\,\lambda_1+\rho,\dots,\lambda_{n-1}+\rho\,\}$. This fact follows
immediately from the following lemma.

\begin{lem}\label{121}
Suppose that the linking matrix of $L$ is 0
and  $\mu_j = \lambda_{k_j} + \rho$, with $j=1,\dots,m$.
Then $ J_L(\mu_1,\dots,\mu_m;q)$ is in 
$q^{(k_{j_1} + \dots + k_{j_m })\frac{n-1}{2}}\Z[q,q^{-1}].$
\end{lem}

\begin{rem} Without the assumption that the linking matrix is 0, the lemma
does not hold true.
\end{rem}

\begin{pf}[Proof of Lemma \ref{121}]
For any framed oriented link $K$ let us defined $J_K^{Hom}(q)$ as follows.
First let  $K'$ be any  link obtained
from $K$ by changing  the framings  so that $K'\cdot K'=0$, i.e.
the sum of all the entries
of the linking matrix of $K'$ is 0. 
Then let $J_K^{Hom}(q)=J_{K'}(\mu_1,\dots,\mu_m;q)$ with all the $\mu_j$ equal
$\lambda_1+\rho$
(i.e. when each $V_{\mu_j}$ is the fundamental representation of $sl_n$).
This is a non-framed version of quantum invariants.
For the link $L$ with 0 linking matrix we choose $L'=L$.
It is known that $\Jhom$ is a version of the Homfly
polynomial and can be calculated by the skein theory as follows.

$$\Jhom_{L_1 \sqcup L_2}(q)= \Jhom_{L_1}(q)\Jhom_{L_2}(q),$$
$$q ^{n/2} \Jhom_{L_+} - q^{-n/2} \Jhom_{L_-} = (q^{1/2} -q^{-1/2})
\Jhom_{L_0},$$
where $(L_0,L_-,L_+)$ is the standard skein triple. 
By using the skein relation and induction on the number of crossing
points of link diagram, one easily proves
\begin{lem} \label{104}If $L$ has $m$ components, then
 $\Jhom_L(q)\in q^{m(n-1)/2} \Z[q^{\pm 1}]$.
\end{lem}
This shows that Lemma \ref{121} holds true when all the $\mu_j$'s
are equal to $\lambda_1+\rho$.

The case when all the $\mu_j$'s are in the set 
$\{\,\lambda_1+\rho,\dots,
\lambda_{n-1}+\rho\,\}$ can be reduced to the case when all $\mu_j$ are equal to
$\lambda_1+\rho$
by using cabling as follows. (See similar arguments in Appendix 1 of
\cite{MasbaumWenzl}).

Let $B_k$ be the braid group on $k$ strands. In what follows we will
consider linear combinations of elements of $B_k$ with coefficients in
$\Z(q^{1/2})$, the ring of rational functions in $q^{1/2}$ with integer coefficients. 
The anti-symmetrizer  $g^{(k)}$ is defined, for $k=1,\dots,n-1$,
by induction as follows
(see \cite{Yokota}):
$g^{(1)}=1\in B_1$,

\begin{equation}\label{anti}
g^{(k)} = \frac{1-q^{-1}}{1-q^{-k}} (g^{(k-1)}\otimes 1)
- q^{(n-1)/2}\frac{1-q^{1-k}}{1-q^{-k}} (g^{(k-1)}\otimes 1)\, 
\sigma_{k-1} \,(g^{(k-1)}\otimes 1),
\end{equation}
where $\sigma_1,\dots,\sigma_{k-1}$ are the standard generators
of the braid group $B_k$, and $z \otimes 1$ is obtained from
$z$ by adding a vertical strand to the right of the braid $z$.

Now suppose $\mu_j= \lambda_{k_j}+\rho$, where $1\le k_j\le n-1$.
Then  (see \cite{Yokota,MasbaumWenzl})
\begin{equation}\label{105}
J_L(\mu_1,\dots,\mu_m;q) = J_{\cal L}(\rho+\lambda_1,\dots,\rho+\lambda_1;q),
\end{equation}
 where
$\cal L$ is a linear combination of links with coefficients in
$\Z(q^{1/2})$ obtained as follows. First we replace the $j$-th
component of $L$ by $k_j$ its push-offs, then we cut these push-offs
at one place and  glue in the  anti-symmetrizer 
element $g^{(k_j)}$. One gets a  linear combinations of links.
Modify the links in this linear combination by changing 
all the framings to 0. The result is $\cal L$. 

Note that since $L$ has 0 linking matrix, the right hand side 
of (\ref{105}) is equal to $\Jhom_{\cal L}(q)$.

{\em Claim:} \qquad \qquad $g^{(k)} = \sum_{z} f_z(q) z,$ \\
where the sum is over a finite subset of the braid group $B_k$, and
$f_z(q)$ is a rational function in
$q^{(1-\sn(z))(n-1)/4} \Z(q)$. Here $\sn(z)$ is
the sign of the permutation corresponding to the braid $z$.
Note that $1-\sn(z)$ is always even.

The claim can be proved easily
by induction on $k$, using (\ref{anti}) and noting that $\sn(zz')=\sn(z)\sn(z')$.

For a braid $z\in B_k$ let $c(z)$ be the number of cycle in 
the permutation corresponding to $z$. The closure of $z$ is a link
of $c(z)$ components.

It is easy to see that, for $z\in B_k$, 
\begin{equation}\label{108}
q^{(1-\sn(z))(n-1)/4} \Z(q)= q^{(k-c(z))(n-1)/2} \Z(q)
\end{equation}

Hence the
right hand side of  (\ref{105}), by Lemma \ref{104}
and  equality (\ref{108}), must belong to
$ q^{(k_{j_1}+\dots +k_{j_m})(n-1)/2}\Z(q)$. On the other hand,
$J_L(\mu_1,\dots,\mu_m;q)$ is a polynomial in $q^{1/2n},q^{-1/2n}$.
It follows that $J_L (\mu_1,\dots,\mu_m;q)$ is in 
$ q^{(k_{j_1}+\dots +k_{j_m})(n-1)/2}\Z[q,q^{-1}]$.
\end{pf}

 \subsection{The Kontsevich integral and weight systems}
The Kontsevich integral is a very powerful invariant of links which
was found by Kontsevich \cite{Kontsevich}. We will use here the 
version for {\em framed } oriented link 
introduced in \cite{LeMurakami2,LeMurakami3}, (with exactly the same
normalization).
In \cite{LeMurakami3} the invariant was denoted by $\hat Z_f(L)$, but 
here for simplicity we will use the notation $Z(L)$. This invariant $Z(L)$
takes values in the (completed) graded vector space
$\cAm$ of Chinese character diagrams
on $m$ circles, i.e. the support of the Chinese character
diagrams is $m$ circles. Here $m$ is the numbers of components of $L$.
See, for example, \cite{LeMurakami3}, for the definition
Chinese character diagrams. In brief, a {\em Chinese character diagram}
consists of a 1-dimensional compact oriented manifold $X$ and
a graph every vertex of which is either univalent or trivalent.
The components of $X$ are supposed be ordered; and
a cyclic order at every trivalent vertex of the graph is fixed.
All the univalent vertex of the graph must be on the manifold $X$.
Usually, the components of $X$ is drawn using solid line, while
the components of the graph -- dashed lines. So we sometime refer to
components of $X$ as solid components, and components of the graph -- dashed.
 A univalent  (resp. trivalent) vertex of the graph is also called an external 
(resp. external) vertex
of the Chinese character diagram. The degree of a Chinese character diagram is
 half the number of  
(internal and external) vertices. The space of
Chinese character diagrams on $X$ is the vector space generated by
Chinese character whose solid part is $X$, subject to the antisymmetry
and Jacobi (also known as IXH) relations.

On the Lie algebra $sl_n$ there is the standard invariant bilinear form
defined by $(y|z)= tr(yz)$, where $tr$ is the trace in the fundamental
representation. 

Suppose $D$ is a Chinese character diagram on $m$ solid circles. Using the 
above defined bilinear form, 
we can define the weight  $W_D(V_1,\dots,V_m)$, which is a number, of the
Chinese character diagram $D$ when the $m$ components of $D$ are colored by
$sl_n$-modules $V_1,\dots,V_m$ (see, for example, 
\cite{Kontsevich,LeMurakami3}).
In fact, $W_D$ can be regarded as a multi-linear map :
$$W_D : \RR^{\otimes m}\to \Q.$$
For $\mu_1,\dots,\mu_m\in \Poo$, let $W_D(\mu_1,\dots,\mu_m)$ stand
for $W_D(V_{\mu_1},\dots,V_{\mu_m})$.

If $Y$ is a linear combination of
chord diagrams on $m$ solid circles, then we can define $W_Y$ using linearity.
The relation between quantum invariants of framed oriented links and the
Kontsevich integral is expressed in the following proposition. 

\begin{pro} \label{Kont}
For $\mu_1,\dots,\mu_m\in \Poo$, one has that 
\begin{equation}
J_L(\mu_1,\dots,\mu_m;q)\bigr|_{q=e^h}= 
\sum _{\ell =0}^\infty W_{Z_\ell(L)}(\mu_1,\dots,\mu_m) h^\ell.
\end{equation}
Here $Z_\ell$ is the degree $\ell$ part of the Kontsevich integral $Z$.
\end{pro}
This is a simple corollary of Drinfeld's theory of quasi-Hopf algebra and was
observed by Le-Murakami \cite{LeMurakami3} and Kassel \cite{Kassel}. 
This shows that all quantum invariants are special values of the Kontsevich integral.

The following is an easy exercise in the theory of the Kontsevich integral
 and we omit the proof (see similar results in \cite{Le1}, Lemma 5.6,
and \cite{Le2}, Proposition 5.3).
\begin{lem}\label{112}
Suppose that $L$ has 0 linking matrix. Then the degree $\ell$ 
part $Z_\ell$ of the Kontsevich integral can be represented in the form
$$ Z_\ell(L)= \sum_D a_D D,  \quad \text{with $a_D\in \Q$} $$
where the sum is over the set of Chinese character diagrams $D$ 
which has at most  $3\ell/2$ external vertices. Moreover, if $s_j$ is
the number of external vertices of $D$ on the $j$-th component, then
$\max_js_j \le \ell -\sum_j s_j$.
\end{lem}
\subsection{Dependence of $W_D$ on $\mu_1,\dots,\mu_m$}
Suppose that $D$ is a Chinese character diagram on $m$  
circles.
At first we will fix $\mu_2,\dots,\mu_m$ and investigate the dependence
of $W_D$ on $\mu_1$, we will write $W_D(\mu_1)$ instead of $W_D(\mu_1,\dots,\mu_m)$.

Break the first component of $D$ at an arbitrary point which is not an external vertex. 
The result  is a Chinese character diagram $D'$ whose support
(i.e. the solid part) is the union of an interval and $m-1$ circle. 
It is known that $D'$, modulo the anti-symmetry and
Jacobi relations, depends only on $D$, but not
on the point where we break the first component (see
 \cite{LeMurakami3}, section 1).

Assign $V_{\mu_2},\dots,V_{\mu_m}$ to the circle components of $D'$,
 and take of the weight of  $D'$. The result is not a number,
 but an element $z(D')$ which lies in the center of the universal enveloping
algebra $U(sl_n)$ of $sl_n$ (see, for example, \cite{Kassel}, Proposition XX.8.2).  
Since $D'$ is determined by $D$, we will write $z(D)$ for $z(D')$.
By the definition of the weight, for every $\mu_1\in \Poo$,
we have that
\begin{equation}\label{110}
W_D(\mu_1)= tr_{\mu_1} z(D),
\end{equation}
where $tr_{\mu_1}$ is the trace taken in the representation 
$V_{\mu_1}$.

Since $z(D)$ is a central element, and $V_{\mu_1}$ is an irreducible
$sl_n$-module, $z(D)$ acts on $V_{\mu_1}$ as  a scalar times the identity. 
The scalar, denoted by $\xi_D(\mu_1)$, is a function on
$\mu_1$, which is known as {\em character} in Lie theory.
Harish-Chandra theory (see, for example, \cite{Humphreys}) says that
$\xi_D(\mu_1)$ is a polynomial function on $\mu_1$ of degree not
exceeding the degree of $z(D)$ in $U(sl_n)$. Moreover the polynomial function
$\xi_D(\mu_1)$ is invariant under  the action of the Weyl group.

\newcommand{\Dim} {\text{\rm dim}}
Formula (\ref{110}) now shows that
\begin{equation} \label{2222}
W_D(\mu_1)= \Dim(V_{\mu_1})\, \xi_D(\mu_1),
\end{equation}
where  the dimension of $V_{\mu_1}$, by the Weyl formula, is equal to $\Delta(\mu_1)$ 
(see Lemma \ref{dimension}).

From the definition of the weight, it follows that 
if $D$ has $k$ external vertices on the first component, then
$z(D)$ has degree $\le k$, and hence $\xi_D(\mu_1)$ is a polynomial
function of degree at most $k$.

We we return to the general case, when the other variables
 $\mu_2,\dots,\mu_m$ are not fixed. Since all the variables
 $\mu_1,\dots,\mu_m$ are independent, we get the following
\begin{pro}\label{111} 
For every $\mu_1,\dots\mu_m\in\Lambda_{++}$, one has that

\begin{equation}
W_D(\mu_1,\dots,\mu_m)=\Delta(\mu_1)\dots \Delta(\mu_m) p_D(\mu_1,\dots,\mu_m),
\end{equation}
where $p_D(\mu_1,\dots,\mu_m)$ is a polynomial function invariant under the
action of the Weyl group.
Moreover, the total degree of $p_D$ is less than or equal to
the number of external vertices of $D$, and the degree of $\mu_j$ in
$p_D$ is less than or equal to the number of external vertices 
on the $j$-th component.
\end{pro}

\subsection{Proof of Proposition \ref{sao22}}

From Proposition \ref{111} and Lemma \ref{112} we see that,
for every $\mu_1,\dots,\mu_m\in\Lambda_{++}$,

$$ W_{Z_\ell(L)}(\mu_1,\dots,\mu_m)= \Delta(\mu_1)\dots\Delta(\mu_m)
p_\ell(\mu_1,\dots,\mu_m),$$
where $p_\ell(\mu_1,\dots,\mu_m)$ is a polynomial function
of degree $\le 3\ell/2$, and if $s_j$ is 
the degree of  $\mu_j$ in $p_\ell$, then $\max_j s_j \le 2\ell -\sum_js_j$.

Applying  Proposition \ref{Kont}, we get

\begin{equation}\label{113}
 J_L(\mu_1,\dots,\mu_m;q)\bigr|_{q=e^h}= \Delta(\mu_1)\dots\Delta(\mu_m)
\sum_{\ell=0}^\infty p_\ell(\mu_1,\dots,\mu_m) h^\ell.
\end{equation}
 
This is true for $\mu_1,\dots\mu_m\in\Lambda_{++}$. We will show that
(\ref{113}) holds true for $\mu_1,\dots\mu_m\in\Lambda$.
 
Suppose  one of the $\mu_j$, say $\mu_1$,  is on the boundary of 
the fundamental chamber $C$.
Then the left hand side is 0 by definition, while
$\Delta(\mu_1)=0$ by Lemma \ref{dimension}. 
Hence (\ref{113}) also holds true in this case.

If we replace $\mu_1$ by $w(\mu_1)$, where $w$ is an element
of the Weyl group,  then the left hand side
of $(\ref{113})$ must be multiplied by $\sn(w)$. 
On the right hand side, all the functions $p_\ell$'s are invariant
under the action of the Weyl group, while
$\Delta(w(\mu_1))= \sn(w) \Delta(\mu_1)$ (see Lemma \ref{dimension}). 
Hence if (\ref{113}) is true 
for $\mu_1$, then it holds true for $w(\mu_1)$. Similarly for other
$\mu_j$. It follows that (\ref{113}) holds true for every
$\mu_1,\dots\mu_m\in \Lambda$.

Now substituting $h= \ln(1+x)$, we get
\begin{align}
Q_L(\mu_1,\dots,\mu_m;q)\bigr|_{q=x+1} &= J_L(\mu_1,\dots,\mu_m;q)\bigr|_{q=x+1}
J_{U^{(m)}}(\mu_1,\dots,\mu_m;q)\bigr|_{q=x+1}\notag \\
&= \Delta^2(\mu_1)\dots\Delta^2(\mu_m)
\sum_{\ell =0}^\infty f'_\ell(\mu_1,\dots,\mu_m) x^\ell,\end{align}
for some polynomial functions $f'_\ell(\mu_1,\dots,\mu_m)$.

Note that $h$, as a formal power series in $x$, has the first
non-trivial term $x$. Hence each $f'_\ell$ is a
polynomial function of total degree $\le 3\ell/2$, and if $s_j$ is the degree
of  $\mu_j$ in $f'_\ell$, then $s_j \le 2\ell -\sum_j s_j$.

That $f_\ell(\mu_1,\dots,\mu_m)\in\Z$ when $\mu_1,\dots\mu_m\in\Lambda$
follows from the fact that $Q_L(\mu_1,\dots,\mu_m)$ is a 
Laurent polynomial in $q$ with {\em integer} coefficients.
This completes the proof of Proposition \ref{sao22}.

\section{Sums related to the root lattice} \label{proof2}
We will always assume that $r$ is an odd prime,
and $\zeta= e ^{2 \pi i/r}$.
Recall that $\rho$ is the half-sum of positive roots, and 
$|\rho|^2 = n(n^2-1)/12$ is always in $\frac{1}{2}\Z$.
If $\beta$ is in the root lattice, then $|\beta|^2$ is an even integer.

\subsection{Gauss sum on the root lattice}
For an integer $b$ not divisible by $r$ let
$$\GG(b)= \sum_{\mu\in(\rho+\Lrootr)} \zeta^{b\frac{|\mu|^2-|\rho|^2}{2}}.$$
Note that for every $\mu\in(\rho+\Lrootr)$, one has that 
$\frac{|\mu|^2-|\rho|^2}{2} \in \Z$, hence $\GG(b)\in\Z[\zeta]$.

The fact that $\zeta$ is an $r$-th of unity has the following 
consequence.
\begin{pro}  \label{210}
For every $\beta$ in the root lattice $\Lroot$, one has that
\begin{equation}
\sum_{\mu\in(\rho+\Lrootr)} \zeta^{b\frac{|\mu+\beta|^2-|\rho|^2}{2} }= \GG(b).
\end{equation}
\end{pro}
\begin{pf}
For $i=1,\dots,n-1$, we have
$$ |\mu+ r \alpha_i|^2-|\rho|^2 = 2r (\mu|\alpha_i) + r^2|\alpha_i|^2.$$
Since $|\alpha_i|^2=2$ is an even number and
$(\mu|\alpha_i)\in\Z$ for every $\mu\in\Lambda$, the right hand side is 
divisible  by $2r$. Hence

\begin{equation}\label{001}
\zeta^{ b\frac{|\mu+r \alpha_i|^2-|\rho|^2}{2} } =
 \zeta^{ b\frac{|\mu|^2-|\rho|^2}{2} }.
\end{equation}
Let  us denote the left hand side of (\ref{210}) by $\GG(b;\beta)$. It's
enough to show that $\GG(b;\beta+\alpha_i)=\GG(b;\beta)$ for every
$\beta\in\Lroot$ and $i=1,\dots,n-1$.
We have 
$$
\sum_{\mu\in(\rho+\Lrootr)} \zeta^{b\frac{|\mu+\beta|^2-|\rho|^2}{2} }=
\sum_{\mu\in(\beta+\rho+\Lrootr)} \zeta^{b\frac{|\mu|^2-|\rho|^2}{2} }
$$
By definition, $\beta+\rho+\Lrootr$ and $\beta+\alpha_i+\rho+\Lrootr$
have the same number of elements, and $\mu$ belongs to the first
set if and only if either $\mu$ or $\mu+ r\alpha_i$ belongs to the
second set (see the definition of
$\Lrootr$ in \S\ref{defi}) Hence, equality (\ref{001}) shows that 
$$\sum_{\mu\in(\beta+\rho+\Lrootr)} \zeta^{b\frac{|\mu|^2-|\rho|^2}{2} }
=
\sum_{\mu\in(\beta+\alpha_i+\rho+\Lrootr)} \zeta^{b\frac{|\mu|^2-|\rho|^2}{2} },$$
or $\GG(b;\beta)=\GG(b;\beta+\alpha_i)$.
\end{pf}

The exact value of $\GG(b)$ is given in the next proposition.
Let us first explain a notation. For a rational number
$a/b$, with $a,b\in\Z$ and $b$ not divisible by $r$, 
let $(a/b)^\vee$ be the natural image of $a/b$ in $\Z/rZ$. 

\begin{pro}
Suppose that $b$ is not divisible by the prime $r$ and $n <r$. Then
$\GG(b)$ has absolute value $r^{(n-1)/2}$, and its phase is
given by
\begin{equation}\label{Gp}
\frac{\GG(b)}{r^{(n-1)/2}}= 
\left( \frac{n}{r}\right)\, 
\left[\left(\frac{b}{r}\right)\, e^{\pi i(1-r)/4}\right]^{n-1}\,
 \zeta^{ (-\frac{ b |\rho|^2}{2})^\vee}.
\end{equation}
Here $\left(\frac{n}{r}\right)$ and $\left(\frac{b}{r}\right)$ are 
the Legendre symbols. In particular,
\begin{equation}
\label{151}
\frac{\GG(b)}{\GG(\sn(b))} = \left(\frac{|b|}{r}\right)^{n-1} \zeta^
{(\frac{\sn(b)-b}{2}|\rho|^2)^\vee}.
\end{equation}
Here $\sn(b)$ is the sign of $b$. 
\end{pro}

\begin{pf}
Recall that $|\rho|^2=n(n^2-1)/12$.
Note that $r+1$ is an even number, and hence $|\rho|^2 (r+1)^2/2$
is an integer. Then we have
\begin{equation} \label{002}
(-\frac{b|\rho|^2}{2})^\vee \equiv -b |\rho|^2 (r+1)^2/2   \pmod{r}.
\end{equation}

Direct calculation shows that
$$\frac{|\mu+\rho|^2-|\rho|^2}{2}  - \frac{|\mu +(r+1)\rho|^2-(r+1)^2
|\rho|^2}{2} = -r (\mu|\rho).$$
For $\mu\in\Lroot$,  $(\mu|\rho)$ is in $\Z$, and hence the right hand side is 
divisible by $r$. It follows that
\begin{align}
\zeta^{b\frac{|\mu+\rho|^2-|\rho|^2}{2} } 
&= \zeta^{b\frac{|\mu+(r+1)\rho|^2-(r+1)^2|\rho|^2}{2}} \notag \\
&= \zeta^{b\frac{|\mu+(r+1)\rho|^2}{2}}\zeta^{(-\frac{b|\rho|^2}{2})^\vee},\label{003}
\end{align} 
where the second identity follows from (\ref{002}). Now we have

\begin{align}
\GG(b) &= \sum_{\mu\in(\rho+\Lrootr)} \zeta^{b\frac{|\mu|^2-|\rho|^2}{2}}
 \notag \\
&= \sum_{\mu\in\Lroot} \zeta^{b\frac{|\mu+\rho|^2-|\rho|^2}{2}} \notag \\
&= \zeta^{(-\frac{b|\rho|^2}{2})^\vee } \sum_{\mu\in\Lroot} 
\zeta^{b\frac{|\mu+(r+1)\rho|^2}{2}}
 \notag \\
&= \zeta^{(-\frac{b|\rho|^2}{2})^\vee } \sum_{\mu\in\Lroot} 
\zeta^{b\frac{|\mu|^2}{2}} \label{005}
\end{align}
where the third identity follows from (\ref{003}) and the fourth from
Proposition \ref{210}, noting that $(r+1)\rho$ is in the root lattice
(since $r+1$ is even).

Recall that $(\alpha_i|\alpha_j)= A_{ij}$, where $A$ is the Cartan matrix.
 If $\mu=k_1\alpha_1+\dots+k_{n-1}\alpha_{n-1}$, then
$|\mu|^2= \bk^t A \bk$, where $\bk$ is the column vector transpose to
$\bk^t= (k_1,\dots,k_{n-1})$. Hence
$$\sum_{\mu\in\Lrootr} \zeta^{b\frac{|\mu|^2}{2}}
= \sum_{\bk \in (\Z/r\Z)^{n-1}} \zeta^{b\, \bk^t A\bk/2}.$$
The right hand side is a multi-variable Gauss sum, and its value is:
\begin{equation}\label{gauss}
\sum_{\bk \in (\Z/r\Z)^{n-1}} \zeta^{b\, \bk^t A\bk/2}=
\left(\frac{n}{r}\right)\, \left[\left(\frac{b}{r}\right)\, e^{\pi i(1-r)/4}\,
 \sqrt r\right]^{n-1}.
\end{equation}
(A proof of this formula is given in Appendix).
This formula, together with (\ref{005}), proves the proposition.
\end{pf}

\begin{cor} \label{1001}
The number $\GG(b)$ is proportional to $\zeta^{(n-1)(r-1)/2}$
by a unit in $\Z[\zeta]$.
\end{cor}
\begin{pf} Note that $\zeta$ is invertible in $\Z[\zeta]$, since
$\zeta^{-1} = \zeta^{r-1}$.  It is well-known that 
$e^{(\pi i(1-r)/4)}\sqrt r$ is proportional to
$\zeta^{(r-1)/2}$ by a unit in $\Z[\zeta]$ (see, for example,
\cite{IK,Murakami}). The corollary now follows immediately from formula (\ref{Gp}).
\end{pf}
\subsection{Completing the square}
\begin{pro} \label{252}
For every $\beta$ in the root lattice $\Lroot$, we have
\begin{equation} \label{222}
\sum_{\mu\in (\rho+\Lrootr)} 
\zeta^{b\frac{|\mu|^2-|\rho|^2}{2} }\zeta^{(\mu|\beta)}
= \zeta^{(-|\beta|^2/2b)^\vee} \GG(b).
\end{equation}
\end{pro}
\begin{pf} The proof uses the trick of completing the
square.
Suppose $b^*$ is a number such that $bb^* \equiv 1  \pmod{r}$. Noting that
 $(\mu|\beta) \equiv b (\mu|b^*\beta)$  $\pmod{r}$, we have
 
\begin{equation}
 \frac{b}{2} (|\mu|^2-|\rho|^2)  + (\mu|\beta) \equiv
\frac{b}{2} (|\mu + b^*\beta|^2 -|\rho|^2) - 
(\frac{|\beta|^2}{2b})^\vee   \pmod{r}.\end{equation}

It follows that the left hand side of (\ref{222}) is equal to
$$
\zeta^{(-|\beta|^2/2b)^\vee} \sum_{\mu\in (\rho+\Lrootr)} 
\zeta^{b\frac{|\mu+b^*\beta|^2-|\rho|^2}{2} }
$$
which is equal to the right hand side of (\ref{222}) by 
Proposition \ref{210}.
\end{pf}

\subsection{The value of $F_{U_b}$}
Recall that $U_b$ is the unknot with framing $b$. We have
$$ J_{U_b}(\mu;q)= q ^{b \frac{|\mu|^2-|\rho|^2}{2}} J_U(\mu;q).$$
Using the value of the unknot (\ref{unknot}) and the definition of $Q_L$ we get
$$ Q_{U_b}(\mu;q)= q ^{b \frac{|\mu|^2-|\rho|^2}{2}} 
q^{2(\mu-\rho|\rho)}\frac{\psi(\mu;q)^2}{\psi(\rho;q)^2}.$$
Hence 
\newcommand{\summuz}{\sum_{\mu\in(\rho+\Lrootr)} \zeta^
{b \frac{|\mu|^2-|\rho|^2}{2}}}
\begin{align} 
F_{U_b}(r)&= \summu Q_L(\mu;\zeta)\notag\\
&=\frac{\zeta^{-2|\rho|^2}}{\psi(\rho;\zeta)^2}
 \summuz [\zeta^{(\mu|\rho)}\psi(\mu;\zeta)] ^2.
\label{011}
\end{align}
\begin{pro} The value of $F_{U_b}(r)$ is given by
\begin{equation}\label{FUp}
F_{U_b}(r)=  n!\, \GG(b) 
\frac{\psi(b^* \rho;\zeta)}{\psi(\rho;\zeta) \psi(-\rho;\zeta)}.
\end{equation}
In particular, when $b=\pm 1$, we have
\begin{equation}\label{Upm}
 F_{U_\pm} = \frac{n!\, \GG(\pm 1)}{\psi(\mp \rho;\zeta)}
\end{equation}
\end{pro}
\begin{pf}

The Weyl denominator formula (see Lemma \ref{weyl}) shows that
\begin{equation}\label{015}
\zeta^{(\mu|\rho)}\psi(\mu;\zeta) = \sum_{w\in\Weyl} \sn(w) \zeta ^{(\mu|w(\rho))}.
\end{equation}
Hence
\begin{equation}\label{010}
[\zeta^{(\mu|\rho)}\psi(\mu;\zeta)]^2 = \sum_{w,w'\in\Weyl} \sn(w w') 
\zeta ^{(\mu|w(\rho)+ w'(\rho))}.
\end{equation}
Although $\rho$ may not be in the root lattice, the sum 
$w(\rho) +w'(\rho)$ is always in $\Lroot$, since each of
$w(\rho),w'(\rho)$ is in $\rho+\Lroot$ (see Lemma \ref{inv}) and 
$2\rho\in\Lroot$.

Using  (\ref{010})  in (\ref{011}) and applying Proposition \ref{252}, we get
\begin{equation}
\label{012}
F_{U_b}(r) = \frac{\zeta^{-2|\rho|^2}\GG(b)}{\psi(\rho;\zeta)^2}
\sum_{w,w'\in\Weyl} \sn(w w') \zeta ^{(-|(w(\rho)+ w'(\rho)|^2/2b)^\vee}.
\end{equation}
Note that 
\begin{align}
|w(\rho)+ w'(\rho)|^2 &= |w(\rho)|^2 + |w'(\rho)|^2 + 2(w(\rho)|w'(\rho))
\notag \\
&= 2 |\rho|^2 + 2(\rho|w^{-1}w'(\rho)). \notag
\end{align}
It's easy to check that, modulo $r$,
$$ \left(\frac{2 |\rho|^2 + 2(\rho| w^{-1}w'(\rho))}{2b}\right)^\vee \equiv 
2 b^*|\rho|^2 + (b^*\rho| w^{-1}w'(\rho)-\rho),
$$
 where $b^*$ is any integer such that 
$b b^* \equiv 1 \pmod r$.
Using this in the right hand side of (\ref{012}), we get
$$
F_{U_b}(r) = \frac{\zeta^{(-2 b^*|\rho|^2 - 2|\rho|^2)}\GG(b)}
{\psi(\rho;\zeta)^2}
\sum_{w,w'\in\Weyl} \sn(w w') \zeta ^{(-b^*\rho| w^{-1}w'(\rho))}.
$$
Note that $\sn(w w')= \sn (w^{-1}w')$. Hence
$$\sum_{w,w'\in\Weyl} \sn(w w') \zeta ^{(-b^*\rho| w^{-1}w'(\rho))}= 
|\Weyl| \sum_{w\in\Weyl} 
\zeta ^{(-b^*\rho| w(\rho))}.
$$
The right hand side, again by Lemma \ref{weyl} and $|\Weyl|= n!$, is equal
to $n!\, \psi(-b^*\rho;\zeta)$.
This shows that 
\begin{equation}\label{1002}
F_{U_b}(r)=  n!\, \GG(b)\, \zeta^{(-2|\rho|^2 - 2 b^*|\rho|^2)}\,
\frac{\psi(-b^* \rho;\zeta)}{\psi(\rho;\zeta)^2}.
\end{equation}
For every number $k$,
\begin{align}
\psi (k \rho;\zeta) &= \prod_{\alpha\in\Phi_+}
(1-\zeta^{-(k\rho|\alpha)}) \notag \\
&= \prod_{\alpha\in\Phi_+}
\zeta^{-(k\rho|\alpha)} \prod_{\alpha\in\Phi_+} 
(\zeta^{(k\rho|\alpha)}-1) \notag \\
&= \zeta^{-2k|\rho|^2} (-1)^{n(n-1)/2}\, \psi(-k\rho;\zeta).\notag
\end{align}
Using this formula, with $k=b^*$ and $k=-1$, in formula (\ref{1002}),
we get (\ref{FUp}).
\end{pf}

\subsection
{Proof of Proposition \ref{sao82}}  
If $k$ is an integer between
$1$ and $r-1$, then $\zeta^k-1$ is  proportional to $\zeta-1$
by a unit in $\Z[\zeta]$. In fact, both 
$$ \frac{\zeta ^k-1}{\zeta-1}  \quad \quad \text{and} \quad \quad
\frac {\zeta-1}{\zeta^k-1} = \frac{\zeta^{k k^*} -1}{\zeta^k-1}$$
are in $\Z[\zeta]$.

For every positive root $\alpha$, the number $(\alpha|\rho)$ is an
integer between 1 and $n-1$. This fact follows from the explicit formulas
of $\alpha$ and $\rho$ (see \cite{Humphreys}).
There are $n(n-1)/2$ positive roots, and hence if $r >n$, the number
$$\psi(\mp\rho;\zeta) = \prod_{\alpha\in\Phi_+}
(1-\zeta^{\pm (\rho|\alpha)}   )$$
is proportional to $(\zeta-1)^{n(n-1)/2}$ in $\Z[\zeta]$.

Formula (\ref{Upm}) and Corollary \ref{1001} shows that
$F_{U_\pm}(r)$ is proportional to 
$$\zeta^{(r-1)(n-1)/2 -n(n-1)/2}$$ 
by a proportional factor which is a unit  in $\Z[\frac{1}{n!}][\zeta]$.
This proves Proposition \ref{sao82}.

\subsection{Simple lens spaces and the constant term}\label{len}
With the knowledge of $F_{U_b}(r)$, we can calculate the value
$\tauPSU_r$ and $\tauPSU$ of lens spaces $M(b)$.

\begin{pro} Let $M(b)$ be the lens space obtained by surgery
on the unknot with framing $b$, where $b$ is a non-zero integer with absolute
value less than the prime $r$. Then
$$\tauPSU_r(M(b))= \left( \frac{|b|}{r}\right)^{n-1} \, 
\zeta^{(\frac{ \sn(b)-b}{2}|\rho|^2)\vee}\,
\frac {\psi(b^*\rho;\zeta)}{\psi(\sn (b) \rho;\zeta)}.$$
\begin{equation}\label{lens}
\tauPSU (M(b))= \left( \frac{|b|}{r}\right)^{n-1} \, 
(1+x)^{(\frac{ \sn(b)-b}{2}|\rho|^2)}\,
\frac {\psi(\rho/b;1+x)}{\psi(\sn (b) \rho;1+x)}.
\end{equation}
\end{pro}
\begin{pf} The first identity follows from 
$\tauPSU_r(M(b))= F_{U_b}(r)/F_{U_{\sn(b)}}(r)$, and
the explicit formula (\ref{FUp}) of $F_{U_b}$. Actually 
the value of $\tauPSU_r(M(b))$ had been calculated by
Takata, see \cite{Takata} for 
the values of $\tauPSU_r$ of Seifert fibered space. 

The value of $\tauPSU (M(b))$ can also be calculated directly from
the explicit formula of $Q_{U_b}(\mu;q)$. However, it's simpler to notice
that the right hand side of (\ref{lens}), which we denote by $g$ for a moment,
satisfies
$${\frak p}_r(g) = {\frak p}_r\left(\tauPSU_r(M(b)\right).$$
Hence by the uniqueness (see Lemma \ref{unique}), we can conclude that
$g = \tauPSU (M(b))$.
\end{pf}

From formula (\ref{lens}) one can calculate the constant term
of the power series $\taut(U_{b})$; the result is 
 $\left( \frac{|b|}{r}\right)^{n-1} \,
|b|^{-n(n-1)/2}$.
Suppose $L'$ is a link with diagonal linking matrix, with $b_1,\dots,b_m$
on the diagonal. From the definition it follows that 
the constant term of the power series $\taut(L')$ is the same
as that of $\taut(U_{b_1}) \times \dots \times \taut(U_{b_m})$.
Hence the constant term of $\taut(L')$ is
is $\left( \frac{|b_1\dots b_m|}{r}\right)^{n-1} \,
|b_1\dots b_m|^{-n(n-1)/2}$. It follows that for every rational
homology 3-sphere $M$, the constant term of the power series 
$\tauPSU (M)$ is 
 $\left( \frac{|H_1(M,\Z)|}{r}\right)^{n-1} \,
|H_1(M,\Z)|^{-n(n-1)/2}$.

\subsection{Proofs of Proposition \ref{new1} and \ref{new2}} 
\label{gammb}
The following follows
directly from the definitions of ${\frak p}_r$.
\begin{lem} \label{new3}
For every $a,b\in\Z$ and every prime $r$ not
 dividing $b$, one has that
$${\frak p}_r\left((1+x)^{a/b}\right) =
{\frak p}_r\left((1+x)^{(a/b)^\vee}\right)= 
{\frak p}_r\left(\zeta^{(a/b)^\vee}\right).$$
\end{lem}

Let $\beta\in\Lroot$ be an element of the root lattice. Then 
Proposition \ref{252} and formula (\ref{Upm}) 
show that
\begin{align}
\frac{\displaystyle{\summu}\zeta^{b\frac{|\mu|^2-|\rho|^2}{2}}
\zeta^{(\mu|\beta)}}
{F_{U_{\sn(b)}}}&=  \frac{1}{n!}\, 
\frac{\GG(b)}
{\GG(\sn(b))} 
\psi\left(-\sn(b)\rho;\zeta\right)\,
\zeta^{(-\frac{|\beta|^2}{2b})^\vee}
\notag \\
&= \left(\frac{|b|}{r}\right)^{n-1}
\frac{\zeta^{(\frac{\sn(b)-b}{2} |\rho|^2)^\vee}}{n!}\,
 \psi\left(-\sn(b)\rho;\zeta\right)\,
\zeta^{(- \frac{|\beta|^2}{2b})^\vee}
\label{1256}
\end{align}
The second equality follows from (\ref{151}).

The right hand side, if the checks are droped and $\zeta$ replaced by
$x+1$, is, by definition, the same as $\Gamma_b(q^{\beta})$. Hence
Lemma \ref{new3} and  (\ref{1256}) show that:
\begin{equation}
\label{sao41}
{\frak p}_r\left(\frac
{\displaystyle{\summu}\zeta^{b\frac{|\mu|^2-|\rho|^2}{2}}
\zeta^{(\mu|\beta)}
}
{F_{U_{\sn(b)}}} \right)= {\frak p}_r(\Gamma_b(q^{\beta})).
\end{equation}
This proves Proposition \ref{new2}, since the set of all
$q^\beta$ spans $\cal S$.

Now we prove Proposition \ref{new1}.
 Recall that 
$$\Gamma_b(q^\beta)= q^{\frac{-|\beta|^2}{2b}}\, y_b,$$
where 
$$ y_b = \left(\frac{|b|}{r}\right)^{n-1}
\frac{(1+x)^{(\sn(b)-b) |\rho|^2/2}}{n!} 
 \, \psi(-\sn(b)\rho;1+x)
$$
does not depend on $\beta$. Since $1-(1+x)^{\sn(b)(\alpha|\rho)}$ is divisible by $x$, and 
since there are $n(n-1)/2$ positive roots, $\psi(-\sn(b)\rho;1+x)$
 is divisible by $x^{n(n-1)/2}$ in $\Z[\frac{1}{n!b}][x]$. Since 
$y_b$ is divisible by $\psi(-\sn(b)\rho;1+x)$, it is divisible by
$x^{n(n-1)/2}$ in $\Z[\frac{1}{n!b}][x]$.

Let $\ab=(a_1,\dots,a_{n-1})$. Then
\begin{align}
{\boldsymbol\eta}^\ab & = \prod_{i=1}^{n-1}(1 -
 q^{-\alpha_i})^{a_i}
\notag \\
&= \prod_{i=1}^{n-1}\sum_{s_i=0}^{a_i} 
{a_i\choose s_i}(-1)^{a_i-s_i}  q^{(s_i-a_i)\alpha_i} 
 \notag \\
&= \sum_{0\le \bs \le \ab} (-1)^{|\ab-\bs|}{\ab\choose\bs} q^{\sum_i 
(s_i-a_i)\alpha_i}. \notag
\end{align}
Here $\bs=(s_1,\dots,s_{n-1})$, $\bs\le \ab$ means $s_i\le a_i$
for every $i$, and ${\ab \choose \bs}$ means $\prod_i {a_i\choose s_i}$.
It follows that
\begin{equation}\label{sao55}
\frac{\Gamma_b(\boldeta^\ab)}{y_b}= (-1)^{|\ab|}\sum_{0\le \bs \le \ab}
 (-1)^{|\bs|} {\ab\choose \bs}
(1+x)^{-|\sum (s_i-a_i)\alpha_i|^2/2b}.
\end{equation}

So we need  only to show that the right hand side of (\ref{sao55})
is divisible by $x^{\lfloor (|\ab|+1)/2\rfloor}$, or that the coefficient of $x^d$ 
in the  right hand side of (\ref{sao55}) is equal to 0
 whenever $2d< |\ab|$. The coefficient of $x^d$ is

$$ v_d =(-1)^{|\ab|}\sum_{0\le \bs \le \ab} (-1)^{|\bs|} {\ab\choose \bs}
{-|\sum (s_i-a_i)\alpha_i|^2/2b \choose d}.$$

The expression
$${-|\sum (s_i-a_i)\alpha_i|^2/2b \choose d} $$
is a polynomial in $s_1,\dots,s_{n-1}$ of total degree $2d$. Hence
when $2d<|\ab|$, $v_d$ is equal to 0 by Lemma \ref{161} below.
This proves Proposition \ref{new1}.

\subsection{A technical result}
The following result is well-known, see \cite{washington}, Lemma 5.19.

\begin{lem}
Suppose that $a$ is a positive integer. Then
$$
\sum_{s=0}^a (-1)^s {a \choose s} s^d=0
$$
for $d=1,2,\dots,a-1$.
\end{lem}
From this lemma one can prove the following.
\begin{lem} \label{161}
Suppose $p(s_1,\dots,s_{n-1})$ is a polynomial
of total degree less than $|\ab|$, where $\ab\in(\Zplus)^{n-1}$. 
Then
$$\sum_{0\le \bs \le \ab} (-1)^{|\bs|}{\ab \choose \bs}\, p(s_1,\dots,s_{n-1})=0.$$
\end{lem}
\begin{pf} It suffices to consider the case when $p(s_1,\dots,s_{n-1})=
s_1^{d_1}\dots s_{n-1}^{d_{n-1}}$. We have

$$\sum_{0\le \bs \le \ab} (-1)^{|\bs|}{\ab \choose \bs} p(s_1,\dots,s_{n-1})=
 \prod_{i=1}^{n-1} \sum_{s_i=0}^{a_i} (-1)^{s_i} {a_i\choose s_i} s_i^{d_i}.$$
Since 
$$d_1 +\dots +d_{n-1} <  a_1+\dots +a_{n-1},$$
there must be an index $i$ for which $d_i<a_i$. Apply the previous
lemma we get the result.
\end{pf}

\subsection{Proof of Proposition \ref{sao71}} Proposition \ref{sao71},
for the case $n=2$, is similar to Lemma 2.2 of \cite{Rozansky}.
We begin with the following lemma.
\begin{lem}\label{991} Suppose $r$ is an odd prime.
For any $d\in\Zplus$, the sum
$ \displaystyle{\sum_{k=0}^{r-1} {k \choose d}} $ is divisible by 
$\displaystyle{(\zeta-1)^{\frac{r-1}{2} -\lfloor \frac{d}{2} \rfloor}}$ 
in $\Zr[\zeta]$, i.e the quotient \quad 
$\displaystyle{\frac{ \sum_{k=0}^{r-1} {k \choose d} }
{(\zeta-1)^{(r-1)/2 -\lfloor \frac{d}{2} \rfloor}}}$ \quad 
is in $\Zr[\zeta]$.
\end{lem}
\begin{pf}
If $d \ge r-1$, then $(r-1)/2 -\lfloor d/2 \rfloor \le 0$. The above
quotient is in $\Z[\zeta]$, and we are done.

Suppose that $d\le r-2$.
Note that  ${ k \choose d}$ is a polynomial in $k$ of degree $d$.
It is well-known that 
\begin{equation}\label{0011}
\sum_{k=0}^{r-1} k^d = \frac{B_{d+1}(r)-B_{d+1}(0)}{d+1},
\end{equation} 
where $B_{d+1}(z)$ is the Bernoulli polynomial (see,
for example \cite{IK}). It is also known that for $d\le r-2$,
the polynomial $B_{d+1}$ has coefficients in $\Zr$. Hence
the right hand side of (\ref{0011}) is divisible by
$r$, which  is divisible by $(\zeta-1)^{(r-1)}$, and hence by
 $(\zeta-1)^{(r-1)/2}$, in
$\Zr[\zeta]$. This proves the lemma.
\end{pf}
\begin{cor}\label{992}
Suppose $p(k_1,\dots,k_{n-1})$ is a polynomial of degree $d$
which takes values in $\Zr$ whenever $k_i\in\Z$. Then
$$\sum_{i=1}^{n-1}\sum_{k_i=0}^{r-1} p(k_1,\dots,k_{n-1})$$
is divisible by 
$\displaystyle{(\zeta-1)^ {(n-1)\frac{r-1}{2} -\lfloor d/2 \rfloor}}$
in $\Zr[\zeta]$.
\end{cor}
\begin{pf}
The polynomial $p(k_1,\dots,k_{n-1})$ is a $\Zr$-linear combination
of terms of the form 
$${k_1\choose d_1} \dots {k_{n-1} \choose d_{n-1}}.$$
Applying Lemma \ref{991}, we get the result.
\end{pf}
[{\em Proof of Proposition \ref{sao71}}].
We have  to show that the number
$$ u = \summu \zeta^{b\frac{|\mu|^2-|\rho|^2}{2}}\, 
{{{\boldsymbol{\alpha}}} \choose \ab}(\mu)$$
is divisible by 
$$(\zeta-1)^{(n-1)\frac{r-1}{2} - \lfloor \frac{|\ab|}{2} \rfloor}$$
in $\Zr[\zeta]$.
Note that
$$u = \sum_{\mu \in\Lrootr} \zeta^{b\frac{|\mu+\rho|^2-|\rho|^2}{2}} 
{{{\boldsymbol{\alpha}}} \choose \ab}(\mu+\rho).$$
Writing $\zeta= 1+(\zeta-1)$ and using the expansion $\zeta^c = \sum_l
{ c \choose l} (\zeta-1)^l$ in the above formula, we can express
$u$ as a polynomial in $\zeta-1$:
$$u =\sum_\ell u_\ell (\zeta-1)^\ell.$$
The coefficient of $(\zeta-1)^\ell$ is
$$u_\ell=
\sum _{\mu\in \Lrootr} {b(|\mu+\rho|^2 -|\rho|^2)/2 \choose \ell}
{{{\boldsymbol{\alpha}}}\choose \ab} (\mu+\rho).$$
Suppose $\mu=k_1\alpha_1 +\dots + k_{n-1}\alpha_{n-1}$, then
$\displaystyle{\sum_{\mu\in\Lrootr}}$ becomes 
$\displaystyle{\sum_{i=1}^{n-1}\sum_{k_i=0}^{r-1}}$, and 
$u_\ell$ can be written as
$$u_\ell = \sum_{i=1}^{n-1}\sum_{k_i=0}^{r-1} p(k_1,\dots,k_{n-1}),$$
where $p(k_1,\dots,k_{n-1})$ is a polynomial taking values
in $\Zr$ whenever $k_i\in\Z$. The degree of $p$ is $2\ell + |\ab|$.
Hence by Corollary \ref{992}, $u_\ell$ is divisible by
$(\zeta-1)^{(n-1)\frac{r-1}{2} - \lfloor \frac{|\ab|}{2} \rfloor-\ell}$ in $\Zr[\zeta]$.

Now remember that $u_\ell$ is the coefficient of $(\zeta-1)^\ell$ in $u$.
Hence taking $(\zeta-1)^\ell$ into account, we see that each term
$u_\ell (\zeta-1)^\ell$, and hence $u$,  is divisible by $(\zeta-1)^{(n-1)\frac{r-1}{2} - 
\lfloor \frac{|\ab|}{2} \rfloor}$ in $\Zr[\zeta]$. 
This completes the proof of Proposition \ref{sao71}.

\appendix
\section {}
\subsection{On the definition of $PSU(n)$-quantum invariant} Here we prove 
that  $\tauPSU_r(M)$ is coincident with the one introduced
by  Kohno and Takata \cite{KohnoTakata2} (see also \cite{Yokota}).
 We first recall Kohno and Takata's definition of quantum $PSU(n)$
invariant, 
which we will denote by $\tilde\tauPSU_r(M)$.

Recall that $\Lambda'_r$ is the intersection of the simplex $\SSS$ and
the weight lattice $\Lambda$.
Let us consider the subset $\Lambda''_r$ of $\Lambda'_r$ which consists of all
weight $\mu=\rho + k_1\lambda_1+\dots +k_{n-1}\lambda_{n-1}$ such that

\begin{equation}
k_1+2k_2+\dots +(n-1)k_{n-1} \quad  \text{is divisible by $n$}\tag{*}
\end{equation} 
Now we define
$F''_L(r)$ by the same formula of $F'_L(r)$, only replacing $\Lambda'_r$ by 
$\Lambda''_r$:

$$F''_L(r)= 
\sum_{j=1}^m\sum_{\mu_j\in\Lambda''_r}Q_L(\mu_1,
\dots, \mu_m;\zeta).$$
Then, by definition (see \cite{KohnoTakata2})

\begin{equation}\label{psln}
\tilde\tauPSU_r(M)= \frac{F''_L(r)}{F''_{U_+}(r)^{\sigma_+}\,
F''_{U_-}(r)^{\sigma_-}}.
\end{equation}

We have the following simple observation.
\begin{lem} The tuple $(k_1,\dots,k_{n-1})$ satisfies $(*)$ if and only if
$k_1\lambda_1+\dots + k_{n-1}\lambda_{n-1}$ is in the root lattice. In other words,
 $$\Lambda''_r=\SSS\cap(\rho+\Lroot).$$
\end{lem}
\begin{pf} If $\beta= l_1\alpha_1+\dots l_{n-1} \alpha_{n-1}$, then
$\beta= k_1\lambda_1+\dots+k_{n-1} \lambda_{n-1}$, where
$\kk = A \bl$. Here $\kk=(k_1,\dots,k_{n-1})$ and $\bl= (l_1,\dots,l_{n-1})$,
 and $A$ is the Cartan matrix.
The inverse  $A^{-1}$ of $A$ is the symmetric matrix
whose entries are given  by $(A^{-1})_{ij} = (n-i)j/n$ for $n-1\ge i\ge j\ge 1$.

Now it is easy to show that $\kk$ satisfies $(*)$ if and only if \/
$\bl= A^{-1}\kk$ has {\em integer} entries.
\end{pf}

The translation group $r\Lroot$ acts on $\Lambda\otimes \R$; and a fundamental
domain is the parallelepiped
$$P_r= \{\, t_1\alpha_1+\dots +t_{n-1}\alpha_{n-1} \mid 
 1\le t_i < r+1, t_i\in\R\,\}.$$
The affine group $\Weylr$ is the semi-direct product of 
the Weyl group $\Weyl$ and
the translation group $r\Lroot$. Hence $r\Lroot$ has 
index $n!$ in $\Weylr$. The group
$\Weylr$ acts on $\Lambda\otimes \R$, and a fundamental domain in $C_r$.

In a sense,
$P_r$ is $n!$ times $C_r$. More precisely, each point in the {\em interior} of $C_r$
has exactly $n!$ points in $P_r$ in its $\Weylr$ orbit
(see \cite{Kac}, Lemma 6.6). For a point on the boundary
of $C_r$, the cardinality of its $\Weylr$ orbits in $P_r$ may be different from
$n!$. However, if one of the $\mu_j$ is on the boundary of
$C_r$, then $Q_L(\mu_1,\dots,\mu_m;\zeta)=0$, by Proposition 
\ref{91}. By that same proposition, $Q_L(\mu_1,\dots,\mu_m,\zeta)$ is invariant
under the action of $\Weylr$.

Notice that $\Lambda''_r= C_r\cap (\rho+\Lroot)$, while $\rho+\Lrootr=
P_r \cap (\rho+\Lroot)$. Hence
$$ \sum_{\mu_j\in (\rho+\Lroot)} Q_L(\mu_1,\dots,\mu_m;\zeta) = (n!)^m 
\sum_{\mu_j\in \Lambda''_r} Q_L(\mu_1,\dots,\mu_m;\zeta).
$$
Or   \quad $ F_L(r) = (n!)^m F''_L(r)$.
When the linking matrix of $L$ has non-zero determinant,
 one has that $\sigma_+ + \sigma_-= m$. Hence
$$\frac{F_L(r)}{F_{U_+}(r)^{\sigma_+}\, F_{U_+}(r)^{\sigma_+}} = 
\frac{F''_L(r)}{F''_{U_+}(r)^{\sigma_+}\, F''_{U_+}(r)^{\sigma_+}}.$$
This shows that our definition  agrees with
that of Kohno and Takata.

\subsection{The value of Gauss sums}
For  the Legendre symbols 
we have $\left(\frac{b}{r}\right)
\left(\frac{b'}{r}\right)\,=\left(\frac{b b'}{r}\right)\,$.
Hence we always have $\left(\frac{b}{r}\right)= 
\left(\frac{b^*}{r}\right)$, where
$b^*$ is any number such that $b b^* \equiv 1 \pmod{r}$.

Let us consider the quadratic Gauss sum
$$\gau = \sum_{k\in \Z/r\Z} \zeta^{k^2},$$
where $\zeta=e^{2\pi i/r}$. The value of $\gau$ is well-known (see,
for example, \cite{IK}): 
$$ \gau=
\begin{cases}
\sqrt r        &  \text{ if $r\equiv 1 \pmod{4}$,}\\ 
 i \,\sqrt r   &  \text{ if $r\equiv 3 \pmod{4}.$}
\end{cases} 
$$
Noting that
$\left(\frac{2}{r}\right) = (-1)^{(n^2-1)/2}$ (see \cite{IK}) we can reformulate
the value of $\gau$ as follows.
\begin{lem} \label{gau}
For an odd prime $r$, the value of the quadratic Gauss sum is given by 
$$\gau = \left(\frac{2}{r}\right) \,  e^{\pi i (1-r)/4}\,  \sqrt r.$$
\end{lem}
The following is well-known (see \cite{IK}, Proposition 6.3.1). 
\begin{lem}\label{gaup}
Suppose that $b$ is an integer not divisible by $r$, then
$$\sum_{k\in\Z/r\Z} \zeta^{b k^2} \, = \left(\frac{b}{r}\right) \gau.$$
\end{lem}
 Now let us consider the following multi-variable Gauss sum
$$\gau_{A,b} = \sum_{\kk \in (\Z/r\Z)^{n-1}}
\zeta^{b \,\kk^t A \kk/2},$$
where $A$ is the Cartan matrix of $sl_n$. Recall that $A_{ii}=2$,
$A_{i,i+1} = A_{i+1,i} = -1$, and other entries are 0.
Since the entries on the diagonal of $A$ are even, 
$\kk^t A \kk/2 $ is always an integer. Hence $\gamma_{A,b}\in\Z[\zeta]$.

Let $D$ be the diagonal $(n-1)\times(n-1)$-matrix with 
entries $D_{jj}=\frac{j+1}{2j}$. Let 
$P$ be the upper-triangle 
$(n-1)\times (n-1)$-matrix with entries $P_{ii}= 1$ for $i=1,\dots,n-1$,
$P_{i,i+1} = - \frac{i}{i+1}$ for $i=1,\dots,n-2$, and other
entries equal to 0. Then it easy to check that
\begin{equation}
\frac{1}{2}A = P^t D P
\end{equation}
i.e., the bilinear form corresponding to $A$ can be diagonalized (over $\Q$)
using $P$. 

Recall that ``check" is the natural reduction modulo $r$ map from
$\Zr$ to $Z/r\Z $.
Since all the entries of $P,D, P^{-1}$ are in $\Zr$ (recall that
$r>n$), there are
defined $D^\vee$ and $P^\vee, (P^{-1})^\vee$. It follows that  the
symmetric bilinear form 
$(A/2)^\vee$ can be diagonalized over $\Z/r\Z$ by $P^\vee$. The resulting
diagonal matrix is $D^\vee$.
So we have
$$ \gau_{A,b} = \sum_{\kk \in (\Z/r\Z)^{n-1}} \zeta ^{ b\, \kk^t D^\vee \kk}.$$
The matrix $D^\vee$ is diagonal, with $D^\vee_{jj} = (j+1) 2 ^* j^*$. Hence
$$ \gau_{A,b}= \prod_{j=1}^{n-1} \sum_{k\in Z/r\Z} \zeta^{b(j+1) 2^* j^* k^2}.$$
Using Lemmas \ref{gau} and \ref{gaup} and the fact that 
$\left(\frac{j}{r}\right)\,\left(\frac{j^*}{r}\right)\,=1$ we get
$$ \gau_{A,b} = \left(\frac{n}{r}\right)\,  
\left[\left(\frac{b}{r}\right)\, e^{\pi i (1-r)/4} \,\sqrt r\right]^{n-1}.$$
This proves formula (\ref{gauss}).
\label{gauss1}

        \end{document}